%% file: article.tex
\documentclass[preprint,sort&compress]{elsarticle}
\usepackage[margin=1in]{geometry}
\usepackage{lineno}
\usepackage{amssymb,amsxtra,amsmath}
\usepackage{latexsym}
\usepackage{ifthen}
\usepackage{graphicx}
\usepackage{calc,fancyhdr}
\usepackage{indentfirst}
\usepackage{parskip,algorithm,algorithmicx}
\usepackage{algpseudocode}
\usepackage{color}
\usepackage{bigints}
\usepackage{subfig}
\usepackage{Article}
\usepackage{color}

\usepackage[usenames,dvipsnames,svgnames,table]{xcolor}
\usepackage{array,multirow,booktabs} 
\usepackage{mathtools}
\newcommand{\R}{\mathbb{R}}
\newcommand{\N}{\mathbb{N}}
\renewcommand{\B}{\mathbb{B}}
\newcommand{\dx}[1]{\mathrm{d} #1}

\DeclareMathOperator*{\argmax}{arg\,max}
\newcommand{\bs}[1]{\boldsymbol{#1}}

\newcommand{\revone}[1]{{\color{black}#1}}
\journal{JCP}

\begin{document}

\begin{frontmatter}

\title{RBF-LOI: Augmenting Radial Basis Functions (RBFs) with Least Orthogonal Interpolation (LOI) for Solving PDEs on Surfaces}

\author[addr1]{Varun Shankar\corref{corresp}}
\address[addr1]{Department of Mathematics, University of Utah, UT, USA}
\ead{vshankar@math.utah.edu}
\cortext[corresp]{Corresponding Author}

\author[addr2]{Akil Narayan}
\address[addr2]{Department of Mathematics, and Scientific Computing and Imaging Institute, University of Utah, UT, USA}
\ead{akil@sci.utah.edu}

\author[addr3]{Robert M. Kirby}
\address[addr3]{Scientific Computing and Imaging Institute, University of Utah, UT, USA}
\ead{kirby@sci.utah.edu}

\begin{abstract}
We present a new method for the solution of PDEs on manifolds $\mathbb{M} \subset \mathbb{R}^d$ of co-dimension one using stable scale-free radial basis function (RBF) interpolation. Our method involves augmenting polyharmonic spline (PHS) RBFs with polynomials to generate RBF-finite difference (RBF-FD) formulas. These polynomial basis elements are obtained using the recently-developed \emph{least orthogonal interpolation} technique (LOI) on each RBF-FD stencil to obtain \emph{local} restrictions of polynomials in $\mathbb{R}^3$ to stencils on $\mathbb{M}$. The resulting RBF-LOI method uses Cartesian coordinates, does not require any intrinsic coordinate systems or projections of points onto tangent planes, and our tests illustrate robustness to stagnation errors. We show that our method produces high orders of convergence for PDEs on the sphere and torus, and present some applications to reaction-diffusion PDEs motivated by biology. 
\end{abstract}

\end{frontmatter}

\input{Intro}
\input{RBFReview}
\input{LOI}
\input{Stability}
\input{Results}
\input{Applications}
\input{Discussion}

\section*{Acknowledgments}

This research was sponsored in part by the Army Research Laboratory and was accomplished under Cooperative Agreement Number W911NF-12-2-0023 . The views and conclusions contained in this document are those of the authors and should not be interpreted as representing the official policies, either expressed or implied, of the Army Research Laboratory or the U.S. Government. The U.S. Government is authorized to reproduce and distribute reprints for Government purposes notwithstanding any copyright notation herein. VS was also partially supported by NSF grants DMS-1521748, CISE AF-1714844, and DMS-1160432, and AN was also partially supported by NSF DMS-1720416 and AFOSR FA9550-15-1-0467.


\section*{References}
\bibliography{article_refs_mod}

\end{document}

%% file: Intro.tex
\section{Introduction}
\label{sec:intro}

Radial Basis Functions (RBFs) are a popular tool for scattered data interpolation and approximation. Much like polynomial approximation methods, RBFs can be used to generate numerical methods for the solution of partial differential equations (PDEs). However, unlike polynomial-based collocation methods, RBF collocation methods are very easily applied to solving PDEs on irregular domains using scattered node layouts~\cite{Bayona2010,Davydov2011,Wright200699}. RBF-based methods also generalize naturally to the solution of PDEs on manifolds $\mathbb{M}\subset\mathbb{R}^3$ using only the Euclidean distance measure in the embedding space and Cartesian coordinates. This ability has been leveraged to obtain four important classes of methods for solving PDEs on manifolds: global RBF methods~\cite{FlyerWright:2007,FlyerWright:2009,FuselierWright:2013,Piret2012}, RBF-Finite Difference (RBF-FD) methods~\cite{FoL11,FlyerLehto2012,SWFKJSC2014,Piret2016}, RBF-Partition of Unity (RBF-PU) methods~\cite{Aiton2011}, and implicit/Hermite RBF-FD methods~\cite{LSWSISC2017}. We will focus on RBF-FD methods for $\mathbb{M} \subset \mathbb{R}^3$ for the remainder of this article.

Historically, the primary drawback of RBF methods has been ill-conditioning in the RBF interpolation matrix~\cite{Fasshauer:2007,Wendland:2004}. This ill-conditioning manifests in algorithmic implementations as a \emph{stagnation} in errors and resulting convergence rates as the number of nodes is increased. In Euclidean domains, this is easily remedied by a change of basis when using Gaussian RBFs~\cite{FornbergPiret:2007,FaMC12,FLF,FoLePo13}, or enforcing reproduction of moderate to high-degree polynomials when using polyharmonic spline (PHS) RBFs~\cite{FlyerPHS,FlyerNS,BarnettPHS,FlyerElliptic}. Unfortunately, none of these methods to offset ill-conditioning appear to apply straightforwardly when the nodes lie on a submanifold $\mathbb{M} \subset \mathbb{R}^d$. Such manifolds appear to require specialized approaches. For instance, Reeger and Fornberg~\cite{ReegerFornbergQuad1,ReegerFornbergQuad2,ReegerFornbergQuad3} are able to overcome these issues on manifolds by using a triangle mesh, projecting the RBF-FD stencil nodes to the tangent plane, and computing the RBF-FD weights there. However, this is a different philosophy from that of using Cartesian coordinates and embedding space distances employed by previous RBF-FD methods on manifolds. We will not pursue this tangent plane approach further in this article. \revone{This choice of working purely on the node set introduces a specific stability issue, which is dealt with in Section 4.2. Yet another way to overcome ill-conditioning with RBF approximations is to view the RBF interpolant as a complex-valued function, and use contour integration~\cite{FoWr} or vector-valued rational approximations~\cite{FoWr2016} in the complex plane to avoid pole singularities. While these alternative techniques may be more applicable to manifolds than the change-of-basis techniques, we leave such an exploration for future work.}

In recent work, however, it was shown that augmenting PHS RBFs with spherical harmonics ameliorated the ill-conditioning issue for interpolation on the sphere $\mathbb{S}^2$~\cite{SWJCP2018}. Recalling that spherical harmonics are merely (orthonormal) restrictions of polynomials in $\mathbb{R}^3$ to $\mathbb{S}^2$, this opens up an alternative approach to overcoming ill-conditioning on manifolds: for an RBF-FD stencil on any manifold $\mathbb{M} \subset \mathbb{R}^3$, compute the (possibly orthonormal) restriction of a polynomial in $\mathbb{R}^3$ to that stencil, and use this polynomial in conjunction with PHS RBFs to enforce polynomial reproduction on that stencil. This is the approach we employ in this article. Of course, this raises the question of how to generate such a restriction. To generate these polynomial subspaces, we turn to Least Orthogonal Interpolation (LOI). LOI is a procedure for generating a minimum-degree polynomial subspace that achieves unisolvency of the interpolation problem associated to a point set in $\R^d$; the point set may have arbitrary finite size and arbitrary geometric configuration \cite{narayan_stochastic_2012}. The LOI procedure requires as input a finite point set and a probability measure on $\R^d$, and outputs a polynomial subspace whose dimension equals the cardinality of the input point set. In addition, LOI also outputs an orthonormal basis for that polynomial subspace.  We note that the LOI procedure is itself a generalization of de Boor and Ron's \emph{least interpolation}~\cite{boor_computational_1992}. Interestingly, RBFs with shape parameters have long been known to recover this least interpolant in a limiting case~\cite{Schaback2005}, further highlighting the connection between LOI and RBF methods. 

We have observed in testing that using LOI polynomial approximations alone in collocation methods for PDEs on scattered nodes produces unstable or unsuitable results. However, this basis has advantageous use in enforcing polynomial reproduction in a PHS-based RBF-FD method, which is the approach we take in this paper. The resulting method, which we call RBF-LOI, retains the strengths of traditional RBF-FD methods on manifolds (Cartesian coordinates and embedding space distance calculations), while overcoming their weaknesses (ill-conditioning and stagnation errors) without the use of higher-precision arithmetic as in~\cite{SWFKJSC2014,LSWSISC2017}. The LOI procedure does increase the cost of the RBF-FD method, as does the growth in stencil sizes required to support RBFs augmented with polynomials. \revone{The cost increases in assembling differentiation matrices corresponding to these larger stencils can be largely ameliorated by using the overlapped RBF-FD method}, a generalization of the RBF-FD method where each local RBF interpolant is used to compute more than one set of RBF-FD weights, thereby drastically decreasing the total number of stencils for a given node set~\cite{ShankarJCP2017}. For this reason, we employ the overlapped RBF-FD method for this article. \revone{We will explore the impact of using overlapped RBF-FD on solution time in Section 5.3. We also remark that while the RBF-LOI procedure appears to be stable on a wide variety of manifolds of different genus, its use on more complicated point-cloud surfaces is a subject of future research.}

The remainder of this paper is organized as follows. In the next section, we present the first application of the overlapped RBF-FD method on manifolds. Section 3 contains a brief description of the LOI procedure for selecting the polynomial subspaces required by the overlapped RBF-FD method. In Section 4, we discuss the stability of the RBF-LOI procedure, and present our approaches to stabilizing our technique. We then validate our methods on the forced heat equation in Section 5 by measuring errors and convergence rates on the sphere and torus. In Section 6, we present applications of our method to solving nonlinear reaction-diffusion equations on more complicated manifolds; more specifically, we solve the Cahn-Hilliard, Fitzhugh-Nagumo, and Turing equations on some interesting manifolds. We conclude with a summary of our results and a discussion of future work in Section 7.

%% file: RBFReview.tex
\begin{table}[h!]
  \begin{center}
  \resizebox{\textwidth}{!}{
    \renewcommand{\tabcolsep}{0.4cm}
    \renewcommand{\arraystretch}{1.3}
    {\scriptsize
    \begin{tabular}{@{}cp{0.8\textwidth}@{}}
      \toprule
      $\bs{x}$ & Point in $\R^d$ \\
      $X$ & Collection of $N$ points $\{\bs{x}_1, \ldots, \bs{x}_N \}$ in $\R^d$\\
      $n$ & Local stencil size \\
      $P_k$ & The collection of the $n$ nearest neighbors in $X$ of $\bs{x}_k \in X$ \\
      $\mathcal{I}_j^k$ & The index in $X$ of the $j$th point in the set $P_k$, with $\mathcal{I}_1^k \equiv k$ \\
      $\mathcal{G}^x$ & $x$-component of the surface gradient\\
      $\rho_k$ & Width of $x_k$-centered stencil $P_k$\\
      $\delta$ & Overlap parameter \\
      $m$ & Polyharmonic spline degree \\
      $R_k$ & Subset of $P_k$ dictated by overlap parameter $\delta$ \\
      $V_P$ & Least orthogonal interpolant polynomial subspace associated to point set $P$ \\
      $M$ & Number of polynomial augmentation terms used in RBF-FD procedure\\
      $h^k_j$ & Ordered orthonormal basis elements for $V_{P_k}$ \\
			$\tau$ & LOI tolerance parameter \\
      $G_k^x$ & RBF-FD weights for the operator $\mathcal{G}^x$ on stencil $k$ \\
    \bottomrule
    \end{tabular}
  }
    \renewcommand{\arraystretch}{1}
    \renewcommand{\tabcolsep}{12pt}
  }
  \end{center}
  \caption{Notation used throughout this article.}\label{tab:notation}
\end{table}

\section{Overlapped RBF-FD on Surfaces}
\label{sec:rbfreview}

\subsection{Description}
\label{sec:description}
We first describe the overlapped RBF-FD method, recently developed by the first author~\cite{ShankarJCP2017}, and its extension to interpolation on $\mathbb{M} \subset \mathbb{R}^3$. Let $X = \{\vx_k\}_{k=1}^N$ be a global set of nodes on the manifold $\mathbb{M}$. Define the stencil $P_k$ to be the set of nodes containing nodes $\vx_{\calI^k_1}$ and its $n-1$ nearest neighbors $\{\calI^k_2,\hdots,\calI^k_n\}$; here, $\{\calI^k_1,\hdots,\calI^k_n\}$ are indices that map into the global node set $X$. We defer discussion of the number of stencils to the end of this section. For the remainder of this discussion, we will focus without loss of generality on the stencil $P_1$. Assume further that we wish to approximate the surface gradient $\nabla_{\mathbb{M}}$, defined in Cartesian coordinates as:
\begin{align}\label{eq:sgrad-op}
\nabla_{\mathbb{M}} = (I - \vn \vn^T)\nabla = [\calG^x,\calG^y,\calG^z]^T,
\end{align}
where $\vn$ is the outward normal, and $\nabla$ is the $\mathbb{R}^3$ gradient. We will first discuss approximating the surface gradient, then use this approximation to approximate the surface Laplacian (Laplace-Beltrami operator),
\begin{align}
\nonumber \Delta_{\mathbb{M}} &= \nabla_{\mathbb{M}} \cdot \nabla_{\mathbb{M}}, \\
\implies \Delta_{\mathbb{M}} &= \calG^x\calG^x + \calG^y\calG^y + \calG^z\calG^z.
\end{align}

Before proceeding, we define the \emph{stencil width} $\rho_1$ as
\begin{align}
\rho_1 = \max\limits_{j=1, \ldots, n} \|\vx_{\calI^1_1} - \vx_{\calI^1_j}\|, 
\end{align}
where $\|\cdot\|$ is the Euclidean norm in $\R^3$. Given an \emph{overlap parameter} $\delta \in (0,1]$, we can now define the \emph{stencil retention distance} $r_1$ to be
\begin{align}
r_1 = (1-\delta)\rho_1.
\end{align}
The parameters $\rho_1$ and $\delta$ define a ball $\mathbb{B}_1$ of radius $r_1$ centered at the node $\vx_{\calI^1_1}$. Let $p_1$ denote the number of nodes in $P_1$ that lie in $\B_1$. Now, let $R_1$ be the set of global indices of the $p_1$ nodes in the subset $\mathbb{B}_1 \subseteq P_1$:
\begin{align}
  R_1 = \{\calR^1_1, \calR^1_2, \hdots, \calR^1_{p_1}\}, 1 \leq p_1 \leq n.
\end{align}
In general, $R_1$ is some subset of a permutation of the global node indices associated to $P_1$. The overlapped RBF-FD method involves computing RBF-FD weights for all the nodes whose indices are in $R_1$, and repeating this process for each stencil $P_k$. Focusing on the $\calG^x$ component of $\nabla_{\mathbb{M}}$, the overlapped RBF-FD weights for all points $\vx \in \mathbb{B}_1$ with indices in $R_1$ are computed using the following augmented local RBF interpolant on $P_1$:
\begin{align}
s_1(\vx,\vy) = \sum\limits_{j=1}^n (g^x)^1_j(\vy) \|\vx - \vx_{\calI^1_j}\|^m + \sum\limits_{i=1}^{M} \lambda^1_i(\vy) h^1_i(\vx),
\label{eq:rbf_interp}
\end{align}
where all superscripts ``1" refer to the stencil index, and $\|\vx - \vx_{\calI^1_j}\|^m$ is the polyharmonic spline (PHS) RBF of degree $m$ ($m$ is odd). A standard RBF-FD procedure would take $h^1_i(\vx)$ as the $M$ monomials corresponding to a carefully chosen $d$-variate polynomial. The key to our technique is the selection of these polynomial basis functions using the least orthogonal interpolation (LOI) technique~\cite{narayan_stochastic_2012}. We will defer discussion on the selection of these basis functions to Section \ref{sec:loi}, and for now consider the $h^1_i$ functions as arbitrary polynomial functions.

The $n$ overlapped RBF-FD weights $(g^x)^1_j(\vy)$ are written explicitly as functions of the evaluation point $\vy$; in other words, each evaluation point $\vy$ in the stencil has a set of $n$ weights associated with it. Our ultimate goal is to compute the weights $L_1$ for the Laplace-Beltrami operator at all points with indices in the set $R_1$. To avoid differentiating normals, we will accomplish this using \emph{iterated interpolation}~\cite{FuselierWright:2013,SWFKJSC2014,LSWSISC2017}. This is done in two steps: first compute overlapped RBF-FD weights for the operators $\calG^x$,$\calG^y$, and $\calG^z$ at all stencil points $\vx_{\calI^1_j}$ (every point in $P_1$); then, combine those RBF-FD weights in such a way that we only compute the weights for all nodes with indices in the set $R_1$. We will now show how these weights are computed for the operator $\calG^x$. We impose the following two (sets of) conditions on the interpolant \eqref{eq:rbf_interp}:
\begin{subequations}
\begin{align}
s_1(\vx_{\calI^1_j},\vx_{\calI^1_i}) &=  \lf.\lf(\calG^x \|\vx - \vx_{\calI^1_j}\|^m\rt)\rt|_{\vx = \vx_{\calI^1_i}}, &i=1,\hdots,n, j=1,\hdots,n, \label{eq:interp_constraint}\\
\sum\limits_{j=1}^n (g^x)^1_j(\vx_{\calI^1_k}) h^1_i(\vx_{\calI^1_j}) &= \lf.\lf(\calG^x h^1_i(\vx)\rt)\rt|_{\vx = \vx_{\calI^1_k}}, & k=1,\hdots,n, i=1,\hdots,M. \label{eq:poly_constraint}
\end{align}
\end{subequations}
The first set of conditions enforces that $s_1(\vx,\vy)$ interpolate the derivatives of the PHS RBF at all the points in $P_1$. The second set of conditions enforces polynomial reproduction/exactness on the overlapped RBF-FD weights. If a degree-$\ell$ polynomial space is employed for $h^1(\vx)$, then $M = {\ell + d \choose d}$; for stability, we also require that $M \leq \lfloor \frac{n}{2} \rfloor$~\cite{ShankarJCP2017,FlyerNS,FlyerPHS}. The interpolant \eqref{eq:rbf_interp} and the two conditions \eqref{eq:interp_constraint}--\eqref{eq:poly_constraint} can be collected into the following block linear system:
\begin{align}
\begin{bmatrix}
A_1 & H_1 \\
H_1^T & O
\end{bmatrix}
\begin{bmatrix}
G^x_1 \\
\Lambda_1
\end{bmatrix}
=
\begin{bmatrix}
B_{A_1} \\
B_{H_1}
\end{bmatrix},
\label{eq:rbf_linsys}
\end{align}
where
\begin{subequations}\label{eq:rbf_linsys_defs}
\begin{align}
  (A_1)_{ij} &= \|\vx_{\calI^1_i} - \vx_{\calI^1_j} \|^m, & i,j&=1,\hdots,n, \\
  (H_1)_{ij} &= h^1_j(\vx_{\calI^1_i}), & i&=1,\hdots,n,\; j=1,\hdots,M,\\
  (B_{A_1})_{ij} &= \lf.\calG^x \|\vx - \vx_{\calI^1_j} \|^m \rt|_{\vx = \vx_{\calI^1_i}}, & i,j& =1,\hdots,n, \\
  \label{eq:h-diff} (B_{H_1})_{ij} &= \lf.\calG^x h^1_j(\vx)\rt|_{\vx = \vx_{\calI^1_j}}, & i&=1,\hdots,M,\; j=1,\hdots,n,\\
  O_{ij} &= 0, & i,j &= 1,\hdots,M.
\end{align}
\end{subequations}
$G^x_1$ is the $n \times n$ local matrix of overlapped RBF-FD weights for the operator $\calG^x$, with each column containing the RBF-FD weights for a point $\vx \in P_1$. The linear system \eqref{eq:rbf_linsys} has a unique solution if the nodes in $P_1$ are distinct~\cite{Fasshauer:2007,Wendland:2004}. More interestingly, \eqref{eq:rbf_linsys} clearly shows that the $n \times M$ matrix of polynomial coefficients $\Lambda_1$ is a set of Lagrange multipliers that enforces the polynomial reproduction constraint \eqref{eq:poly_constraint}. The above procedure can be repeated with the operators $\calG^y$ and $\calG^z$ to obtain the local differentiation matrices $G^y_1$ and $G^z_1$. Next, define the \emph{truncated} matrix $\tilde{G^x_1}$ as:
\begin{align}
\lf(\tilde{G^x_1}\rt)_{ij} = \lf(G^x_1\rt)_{ij}, i=1,\hdots,n, j=1,\hdots,p_1,
\end{align}
\emph{i.e.}, the $n \times p_1$ submatrix of $G^x_1$ corresponding to the nodes in the ball $\mathbb{B}_1$. Similarly define the truncated matrices $\tilde{G^y_1}$ and $\tilde{G^z_1}$. Finally, we use iterated interpolation to obtain the differentiation matrix $L_1$ for the Laplace-Beltrami operator $\Delta_{\mathbb{M}}$ as:
\begin{align}
L_1 = \lf(\tilde{G^x_1}\rt)^T \lf(G^x_1\rt)^T + \lf(\tilde{G^y_1}\rt)^T \lf(G^y_1\rt)^T + \lf(\tilde{G^z_1}\rt)^T \lf(G^z_1\rt)^T.
\end{align}
This construction using truncated matrices ensures that the $p_1 \times n$ matrix $L_1$ only contains RBF-FD weights for the nodes in $P_1$ whose indices are in the set $R_1$. By construction, the rows of $L_1$ populate the rows of a global differentiation matrix $L$, while the \emph{columns} of $\tilde{G^x_1}$ and its counterparts populate the rows of the global differentiation matrices $G^x$, $G^y$, and $G^z$. If the weights for the Laplace-Beltrami operator are not required, it is straightforward to directly compute the truncated matrices $\tilde{G^x_1}$, $\tilde{G^y_1}$, and $\tilde{G^y_1}$ by modifying \eqref{eq:interp_constraint}--\eqref{eq:poly_constraint}.

To avoid computing multiple sets of RBF-FD weights for a node $\vx_k$, we also require that weights computed for some node $\vx_k$ never be recomputed by some other stencil $P_i, i\neq k$. The entire procedure above must be performed for each stencil; this can be computationally onerous if the number of stencils is comparable to the total number of points $N$. Denote the total number of stencils by $N_{\delta}$. For a quasi-uniform node set, $N_{\delta} = \frac{N}{p}$, where $p = \max \lf((1-\delta)^d n,1 \rt)$, and $d$ is the dimension (in the above discussion, $d = 3$). If $\delta=1$, this gives us $N_{\delta} = N$, recovering the standard RBF-FD method. However, if $\delta <1$, then $N_{\delta} << N$, giving a significant speedup over the standard RBF-FD method. For a detailed complexity analysis, see~\cite{ShankarJCP2017}.

\subsection{Parameter selection}
\label{sec:param}

In this section, we describe parameter selection for our method. Given a linear operator $\calL$ of order $\theta$ and an RBF-FD differentiation rule $(\bs{x}_j, w_j)_{j=1}^n$, we have the following error estimate for RBF-FD based formulas that reproduce a polynomial of degree $\ell$~\cite{DavydovSchaback2017}:
\begin{align}
|\calL f(\vx) - \sum_{j=1}^n w_j f(\vx_j)| \leq C(m,\vx) h^{\ell+1-\theta},
\end{align}
where $h$ is the fill distance of the node set, and $m$ is the degree of the RBF PHS used. While a derivation of such a formula for RBF-FD on manifolds is pending, we nevertheless use this formula to guide our parameter selection. If we require an RBF-FD method with order of accuracy $\xi$, we set
\begin{align}
\ell = \xi + \theta - 1.
\label{eq:select_ell}
\end{align}
It is important to note that in the context of the LOI procedure (Section 3), the input value of $\ell$ is merely a ``requested'' polynomial degree. In practice, the LOI procedure may output a polynomial of degree slightly lower than $\ell$. Thus, while the number of polynomial basis functions is related to the input degree $\ell$ as $M = {\ell + d \choose d}$, the LOI procedure (and its accompanying stabilization techniques) may in practice result in a smaller $M$ than requested. Regardless, since this only affects the polynomial reproduction, we select the stencil size as
\begin{align}
n = 2M + 1 =  2 {\ell + d \choose d} + 1.
\label{eq:select_n}
\end{align}
The degree $m$ of the PHS RBF can either be fixed~\cite{FlyerPHS,BarnettPHS,ReegerFornbergQuad1,ReegerFornbergQuad2,ReegerFornbergQuad3}, or increased with respect to $\ell$~\cite{SWJCP2018,iske2002}. In Euclidean domains, it appears beneficial to fix $m$~\cite{ShankarJCP2017,FlyerPHS,FlyerNS,BarnettPHS}. On the other hand, the traditional scaling law $m = 2\ell+1$, appears to give the greatest accuracy and stability on manifolds~\cite{SWJCP2018}. We have found that $m=2\ell+1$ was the most stable choice for all manifolds considered in this article.

Finally, we must also choose the overlap parameter $\delta \in (0,1]$. In practice, we have observed that setting $\delta \leq 0.2$ typically completely decouples the stencils, resulting in ill-posed subproblems when solving PDEs. However, using $0.2< \delta \leq 1$ appears to be perfectly stable, and $\delta$ can be chosen to be smaller for larger values of $n$~\cite{ShankarJCP2017}. Given these constraints, we use the following heuristic:
\[
 \delta =
  \begin{cases} 
      \hfill 0.7 \hfill & \text{ if $\ell \leq 4$} \\
      \hfill 0.5 \hfill & \text{ if $4 < \ell \leq 6$} \\
			\hfill 0.3 \hfill & \text{ if $\ell > 6$} \\
  \end{cases}
\]
We find that these values of $\delta$ result in stable differentiation matrices, while also facilitating the rapid assembly of these matrices.

%% file: LOI.tex
\section{Least Orthogonal Interpolation}
\label{sec:loi}

The polynomially-augmented RBF-FD procedure described above requires specification of the polynomial functions portion of the algorithm, i.e., specification of the functions $h^k_j(\vx)$, $j=1, \ldots, M$. We define this polynomial basis in this section; for simplicity we omit all notational dependence on the stencil index $k$ in this section.

As mentioned in Section \ref{sec:intro}, the LOI procedure outputs a polynomial subspace for a given finite input point set. In addition, the procedure also outputs a basis $h_j(\vx)$ for the polynomial subspace whose elements are orthonormal in a weighted $L^2$ space on $\R^d$, where the weight is given by the differential of a user-prescribed probability measure. This is the basis we will use to compute augmented RBF-FD weights in \eqref{eq:rbf_linsys}. We start our discussion by assuming that the probability measure is given and fixed (denoted $\mu$ below), and describe towards the end of this section our choice for this measure.

\subsection{Notation}
With $d \in \N$, a point $\bs{x} \in \R^d$ has Cartesian components
\begin{align*}
  \vx = \left\{ x^{(1)}, \ldots, x^{(d)} \right\}.
\end{align*}
We use standard multi-index notation: given a multi-index $\alpha \in \N_0^d$, we have
\begin{align*}
  \alpha &= \left( \alpha^{(1)}, \ldots, \alpha^{(d)} \right) \in \N_0^d,  &
  |\alpha| &= \sum_{q=1}^d \alpha^{(q)} &
  \vx^\alpha &= \prod_{q=1}^d \left( x^{(q)}\right)^{\alpha^{(q)}}.
\end{align*}
We use $V_{n-1}$ to denote the space of polynomials of degree $n-1$ or less in $\R^d$:
\begin{align*}
V_{n-1} &= \mathrm{span} \left\{ \vx^\alpha \;\; \big| \;\; \alpha \in \N_0^d, \;|\alpha| \leq n-1 \right\}, & \dim V_{n-1} &= \left(\begin{array}{c} n-1+d \\ d \end{array}\right).
\end{align*}
Let $\mu$ be a probability measure on $\R^d$, and let $L^2_\mu(\R^d)$ be the space of real-valued square-integrable functions with respect to the measure $\mu$ on $\R^d$.  We assume that $\mu$ has finite polynomial moments of all orders and is not degenerate with respect to polynomials, i.e., 
\begin{align}\label{eq:mu-assumptions}
  0 &< \int_{\R^d} \vx^{2 \alpha} \dx{\mu}(x) < \infty, & \alpha & \in \N_0^d
\end{align}
Formally, we require only finite moments up to a finite order for the procedure we discuss, and degeneracy is neither a mathematical nor a computational issue. However, the stronger assumptions above are sufficiently general for our presentation, and codifying the allowable relaxation of the above requirements involves unnecessary technical discussions.

Under the above conditions, there is a sequence of polynomials $\phi_\alpha$, $\alpha \in \N_0^d$, with $\deg \phi_\alpha = |\alpha|$, satisfying
\begin{align*}
  \int_{\R^d} \phi_\alpha(\vx) \phi_\beta(\vx) \dx{\mu}(\vx) &= \delta_{\alpha, \beta}, & V_{n} = \mathrm{span} \left\{ \phi_\alpha \;\; \big| \;\; \alpha \in \N_0^d, \; |\alpha| \leq n \right\}.
\end{align*}
Note that $\mu$ may have compact support, in which case all integrals can be reduced to ones over this compact set. For each $n$, the space $V_n$ is a finite-dimensional Hilbert space.

Assuming polynomials are complete in $L^2_\mu$, any $f \in L^2_\mu$ has the Fourier-like expansion
\begin{align*}
  f(\vx) &= \sum_{\alpha \in \N_0^d} \widetilde{f}_\alpha \phi_{\alpha}(\vx) = \sum_{j=0}^\infty \sum_{|\alpha| = j} \widetilde{f}_\alpha \phi_{\alpha}(\vx) \eqqcolon \sum_{j=0}^\infty f_j(\vx), &
  \widetilde{f}_\alpha &= \int_{\R^d} f(\vx) \phi_\alpha(\vx) \dx{\mu}(\vx),
\end{align*}
where we have defined $f_j \in V_{j}$ in terms of the coefficients $\widetilde{f}_\alpha$. We define the operation $\left(\cdot\right)_{\downarrow}$ as follows:
\begin{align*}
  f_{\downarrow} &\coloneqq f_r(\vx), & r &= \min \left\{ j \in \N_0 \;\; \big| \;\; f_j \neq 0 \right\}.
\end{align*}
Note that this operation depends on $\mu$.

\subsection{The least orthogonal interpolant}
Let $P = \left\{\vx_1, \ldots, \vx_n\right\} \subset \R^d$ be a point set of size $n \in \N$. With $\mu$ fixed, the least orthogonal interpolation procedure provides a polynomial subspace of dimension $n$ associated to $P$.

The Riesz representor of the point-evaluation map $v \mapsto \delta_{\vx_j}(v) = v(\vx_j)$ in the finite-dimensional space $V_{n-1}$ has the form
\begin{align}\label{eq:vj-def}
  v_j(\vx) = \sum_{|\alpha| \leq n-1} \phi_\alpha(\vx_j) \phi_\alpha(\vx).
\end{align}
With $v_j$, $j=1, \ldots, n$, defined above in terms of the nodes in $P$, we can define the following space of polynomials:
\begin{align}\label{eq:LOI-space}
  V_P \coloneqq \mathrm{span} \left\{ v_{\downarrow} \;\; \big|\;\; v \in \mathrm{span}\left\{v_1, \ldots, v_n\right\}\right\}.
\end{align}
The main result from \cite{narayan_stochastic_2012} is that the space $V_P$ has dimension $n$ and the interpolation problem on $P$ in the space $V_P$ is unisolvent. Therefore, we can always identify a unique polynomial in $V_P$ given data on $P$; this polynomial is the least orthogonal interpolant, and $V_P$ is the least orthogonal interpolant (polynomial) space associated to $P$. Because $V_P$ has dimension $n$, there is an $L^2_\mu$-orthonormal basis, $h_1(\cdot), \ldots, h_n(\cdot)$, for $V_P$. This basis can be computationally generated using linear algebra, and these are the basis elements that we use in the RBF-FD procedure \eqref{eq:rbf_linsys} and \eqref{eq:rbf_linsys_defs}. This algorithmic construction depends on detecting rank-deficiency of certain submatrices of a Vandermonde-like matrix. Like all numerical linear algebraic methods to detect rank, this in turn depends on user specification of a tunable tolerance parameter denoted $\tau$. This tolerance parameter is related to numerical unisolvency of the interpolation problem, and we discuss it in more detail in section \ref{sec:stab-tau}.

The LOI definition above is abstract but, given a point set $P$ and data on $P$, the computation of the interpolant (and the basis $h_j$) involves only standard tools from numerical linear algebra, namely $L U$ and $Q R$ factorizations~\cite{narayan_stochastic_2012}. In particular, the computational complexity of the procedure is comparable to that for a standard interpolation problem of size $n$.

\subsection{Differentiation of the interpolant}
Condition \eqref{eq:h-diff} in the formation of local RBF-FD weights shows that we must have the ability to differentiate the polynomial basis $h_j$. The least orthogonal interpolant associated to the points set $P$ of size $n$ has the form
\begin{align*}
  p(\vx) &= \sum_{j=1}^n c_j h_j(\vx) \in V_P, & \int_{\R^d} h_j(\vx) h_\ell(\vx) \dx{\mu}(\vx) = \delta_{j,\ell}
\end{align*}
The orthonormal basis $h_j$ is output from the LOI procedure. (The coefficients $c_j$ may also be easily generated if data on $P$ is given.) The basis elements are polynomials, each of degree $\deg h_j$, and have an expansion in terms of the original orthonormal basis $\phi_\alpha$:
\begin{align*}
  h_j(\vx) = \sum_{|\alpha| = \deg h_j} d_{j,\alpha} \phi_\alpha(\vx),
\end{align*}
for some constants $d_{j,\alpha}$ that are output from the LOI procedure. (Note above that we only need take $\alpha$ satisfying $|\alpha| = \deg h_j$, not $|\alpha| \leq \deg h_j$.) For the purposes of solving PDEs, we are particularly interested in differentiating the interpolant $p(x)$, say with respect to coordinate $\ell$. This is given by
\begin{align*}
\frac{\partial p}{\partial x^{(\ell)}} = \sum_{j=1}^n c_j \sum_{|\alpha| = \deg h_j} d_{j,\alpha} \frac{\partial \phi_\alpha}{\partial x^{(\ell)}}.
\end{align*}
Therefore in order to compute derivatives, we need only the ability to construct the interpolant and to differentiate the original basis $\phi_\alpha$. The Cartesian derivatives of the basis $\phi_\alpha$ can then be combined appropriately to give the surface gradient $\nabla_{\mathbb{M}} \phi_\alpha$.

\subsection{The measure $\mu$}
The LOI procedure is a well-defined map from a point configuration to the sought basis $h_j$ (and its derivatives). However, we have yet to make a specification for the probability measure $\mu$. In principle any measure satisfying \eqref{eq:mu-assumptions} will suffice, but in our quest for stable methods, it seems more prudent to choose $\mu$ so that the input orthonormal basis $\phi_\alpha$ is well-behaved. 

For simplicity, we choose $\mu$ to be the tensor-product Chebyshev measure over the smallest bounding box for the nodal set $P$. With this, we can generate the basis $\phi_\alpha$ as 
\begin{align*}
  \phi_\alpha(\vx) = \prod_{j=1}^d T_{\alpha^{(j)}}\left( x^{(j)}\right),
\end{align*}
where $\{T_q\}_{q=0}^\infty$ are the univariate orthonormal Chebyshev polynomials. In this way, the basis $\phi_\alpha$ along with its partial derivatives are easily computed.

In order to apply all of the above to the RBF problem, for each stencil (nodal set) $P_k$ we compute LOI basis functions $h^k_j$ for use on that stencil, using the prescription of $\mu$ above.

%% file: Stability.tex
\section{Eigenvalue Stability}
\label{sec:stab}

Eigenvalue stability is essential for stable time integration of PDEs. In this context, a discrete version of a diffusive (elliptic) differential operator can be regarded as ``stable'' if its spectrum contains no eigenvalue with positive real parts. In~\cite{SWFKJSC2014}, this stability was achieved by performing a nonlinear optimization for the shape parameter on each stencil, constrained so that the RBF interpolation matrices on each stencil should have approximately the same condition number. In~\cite{LSWSISC2017}, this stability was achieved by encouraging diagonal dominance in RBF-HFD differentiation matrices using a stencil selection algorithm. 

In contrast, stability in the RBF-LOI procedure is primarily achieved by picking a single tolerance parameter $\tau$ that is globally defined over the whole mesh. We also utilize an additional empirical stabilization procedure based on avoiding ``axis alignment", that is alignment of a stencil configuration with the global Cartesian coordinate system. This last correction is a pathology of our Cartesian representation of points on a manifold, and is needed only on stencils for which one or more points geometrically aligns with a Cartesian axis. We have observed that this stabilization is needed on very few stencils in all our tests.

We will now describe both techniques for stabilization. Throughout this section we use the same notation as in Section \ref{sec:loi}: $P$ is a generic RBF-FD stencil containing the points $\vx_1, \ldots, \vx_n$. The principles discussed here are applied for each local RBF-FD stencil.

\subsection{The LOI tolerance parameter}\label{sec:stab-tau}
The tolerance parameter $\tau$ that we use is embedded in the LOI construction procedure. This tolerance parameter is a way of tuning selection of the polynomial degree so that the polynomial interpolation on $P$ is numerically stable.

The mathematical LOI procedure described in Section \ref{sec:loi} can be implemented as a sequence of $Q R$ decompositions \cite{narayan_stochastic_2012}, and the complexity of the entire procedure is asymptotically the same as a standard, square, interpolation problem of the same size, \emph{i.e.}, $O(n^3)$. Recall that, given a finite point set $P$, the LOI polynomial space $V_P$ has dimension $n = |P|$, and the elements $\left\{ h_j \right\}_{j=1}^n$ are an orthonormal basis for $V_P$. 


The algorithmic implementation of the LOI procedure builds the functions $h_1, \ldots, h_k$ sequentially by identifying them with points in $P$, \emph{i.e.}, $h_1, \ldots, h_k$ are built and are each identified with a sequence $Y_k \coloneqq \{\vx_{i_1}, \ldots, \vx_{i_k} \} \subset P$, where $i_1, \ldots, i_k \in \{ 1, \ldots, n \}$. After $h_1, \ldots, h_k$ have been built, then $h_{k+1}$ is identified and constructed with an associated point $\vx_{i_{k+1}} \in P \backslash Y_k$. The identification and construction of $h_{k+1}$ is based on the residual of a projection. Let $k < n$, and define $r \coloneqq \deg h_k$; we must have that $r \leq n-1$. We require truncations of the summation in \eqref{eq:vj-def}:
\begin{align*}
  v_{j,r}(\vx) &= \sum_{|\alpha| \leq r} \phi_\alpha(\vx_j) \phi_\alpha(\vx), & r &\leq n-1.
\end{align*}
Now define $\Pi^\perp_{r}$ as the $L^2_\mu$-orthogonal projector onto orthogonal complement of $V_r$, and let $I_k$ denote the LOI interpolation operator associated with $Y_k$ with the basis $h_1, \ldots, h_k$; i.e., 
\begin{align*}
  I_k f(\vx) &= \sum_{j=1}^k c_j h_j(\vx), & I_k f(\vx_j) &= f(\vx_j), \hskip 5pt j = 1, \ldots, k.
\end{align*}
We then define the degree $r$ ``residual" as
\begin{align}\label{eq:R-def}
  R = \max_{j \in \{1, \ldots, n\}\backslash \{i_1, \ldots, i_k \}} \left\| \Pi^\perp_r \left( v_{j,r}  - I_k v_{j,r} \right) \right\|_{L^2_\mu}.
\end{align}
Now let $\tau > 0$ be a tolerance parameter. This tolerance paramter is a threshold for the allowable residual value $R$. If $R \geq \tau$, then we choose
\begin{align*}
  i_{k+1} = \argmax_{j \in \{1, \ldots, n\}\backslash \{i_1, \ldots, i_k \}} \left\| \Pi^\perp_r \left( v_{j,r}  - I_k v_{j,r} \right)\right\|_{L^2_\mu},
\end{align*}
and $h_{k+1}$ is chosen as a normalized version of $\Pi^\perp_r (v_{i_{k+1},r} - I_k v_{i_{k+1,r}})$.  Otherwise, if $R < \tau$, we set $r \gets r + 1$, recompute $R$ from \eqref{eq:R-def}, and repeat the comparison of $R$ to $\tau$. When $\tau = 0$, only a pathological prescription of $\mu$ allows more than one increment of $r$ for each $k$.

The tolerance parameter $\tau$ can now be understood in terms of $R$. The quantity $R$ measures the ability of the point set $P$ to resolve (with respect to the measure $\mu$) a certain subspace of polynomials of degree $r$. The comparison of $R$ with $\tau$ then enforces a desired threshold of resolvability for subspaces included in the LOI procedure. When this threshold is not met, the LOI interpolation operator will be (relatively) ill-conditioned on this subspace of polynomials. Therefore, instead of including this subspace, we simply increment the degree (generate a new subspace) in hopes of achieving a more stable polynomial approximation. This simple heuristic allows us, for a fixed $\tau$, to achieve stable RBF-LOI approximations for general test cases without any optimization.

\subsection{Axis misalignment}
We empirically observe that some local stencils produce unstable results when the point set $P$ has a very special configuration in space. This instability is not directly caused by our procedures, but instead by our choice of the alignment of a Cartesian coordinate system in $d$-dimensional space; \revone{it is plausible that this instability may not occur, say, for approximations on the tangent plane}. We first describe the source of the instability, and then describe our simple computational strategy to circumvent the issue. To keep notational jargon at a minimum in this section, with $d=3$ we use the notation $\left( x^{(1)}, x^{(2)}, x^{(3)} \right) = \left( x, y, z \right)$.

The stencil $P$ contains a spatial configuration of points, and the LOI procedure outputs a basis $h_1, \ldots, h_M$ from these points. The instability we observe stems from situations where a special arrangement of points $P$ results entries of the matrix defined in \eqref{eq:h-diff} satisfying $\mathcal{G}^x h_j \equiv 0$. This ``zero column" of the matrix $B_{H_1}$ causes numerical instabilities when the corresponding global (sparse) differentiation matrix is used for the discretization of PDEs. 

We give a brief account of why this ``zero column" occurs: The functions $h_j$ are arranged in order of increasing polynomial degree. E.g., in $d=3$, $h_1$ is the constant polynomial, and $h_2$, $h_3$, and $h_4$ are all linear polynomials, except in pathological arrangements of $P$ or for pathological measures $\mu$. For simplicity of discussion, we assume in this section that $d=3$ and that $h_2$, $h_3$, and $h_4$ are all linear polynomials.

Suppose that the outward pointing normal vector $\vn$ at the stencil center equals $(n_1, n_2, n_3)^T$. This implies that the $x$-component $\mathcal{G}^x$ of the surface gradient operator in \eqref{eq:sgrad-op} is given by
\begin{align*}
\left(\begin{array}{c} \mathcal{G}^x \\ \mathcal{G}^y \\ \mathcal{G}^z \end{array}\right) = \left(\bs{I} - \vn \vn^T \right) \nabla
\end{align*}
A linear polynomial, say $h_2$, has the expansion
\begin{align*}
h_2(\vx) = \alpha x + \beta y + \gamma z + d,
\end{align*}
where $\alpha, \beta, \gamma$, and $d$ are all constants. Then we can compute
\begin{align*}
  \left(\begin{array}{c} \mathcal{G}^x \\ \mathcal{G}^y \\ \mathcal{G}^z \end{array}\right) h_2 = \left( \bs{I} - \vn \vn^T \right) \left( \begin{array}{c} \alpha \\ \beta \\ \gamma \end{array}\right).
\end{align*}
We can see then that one component of this vanishes when any row of $\left(\bs{I} - \vn \vn^T \right)$ is orthgonal to $\left( \alpha, \beta, \gamma\right)^T$. While this situation happens rarely, it is not difficult to construct situations when such a condition is triggered. Indeed, by constructing a stencil arranged with a normal vector $\vn = (1, 0, 0)^T$, then certain configurations of the stencil $P$ cause $(\alpha, \beta ,\gamma)^T = (\alpha, 0, 0)^T$, which then results in $\mathcal{G}^x h_2 \equiv 0$. 

We observe in practice that stencils centered on points that are ``axis aligned" with the Cartesian grid generate an LOI basis that causes one surface gradient component for linear polynomials to vanish. Our simple fix to circumvent this alignment issue for those stencils is as follows:
\begin{enumerate}
\item Check if any column of the LOI $\mathbb{R}^3$ gradient-component matrices on the stencil $P_k$ (other than the first) contains only zeros. Mark the columns.

\item Eliminate the basis functions $h^k_j(\vx)$ (and their derivatives) corresponding to these marked columns from the matrices $H_k$ and $B_{H_k}$.
\end{enumerate}

%% file: Results.tex
\section{Results}
\label{sec:results}

We have completed a full description of the RBF-LOI procedure and in this section we test the convergence rates of this method on the sphere and torus, where explicit expressions for surface differential operators are known. We then present timing results for our methods as a function of error. In all cases, the relative errors are measured in the embedding space and no quadrature is used. This may affect the constants involved, but should not affect the convergence rates. 
\begin{figure}[h!]
\centering
\includegraphics[scale=0.6]{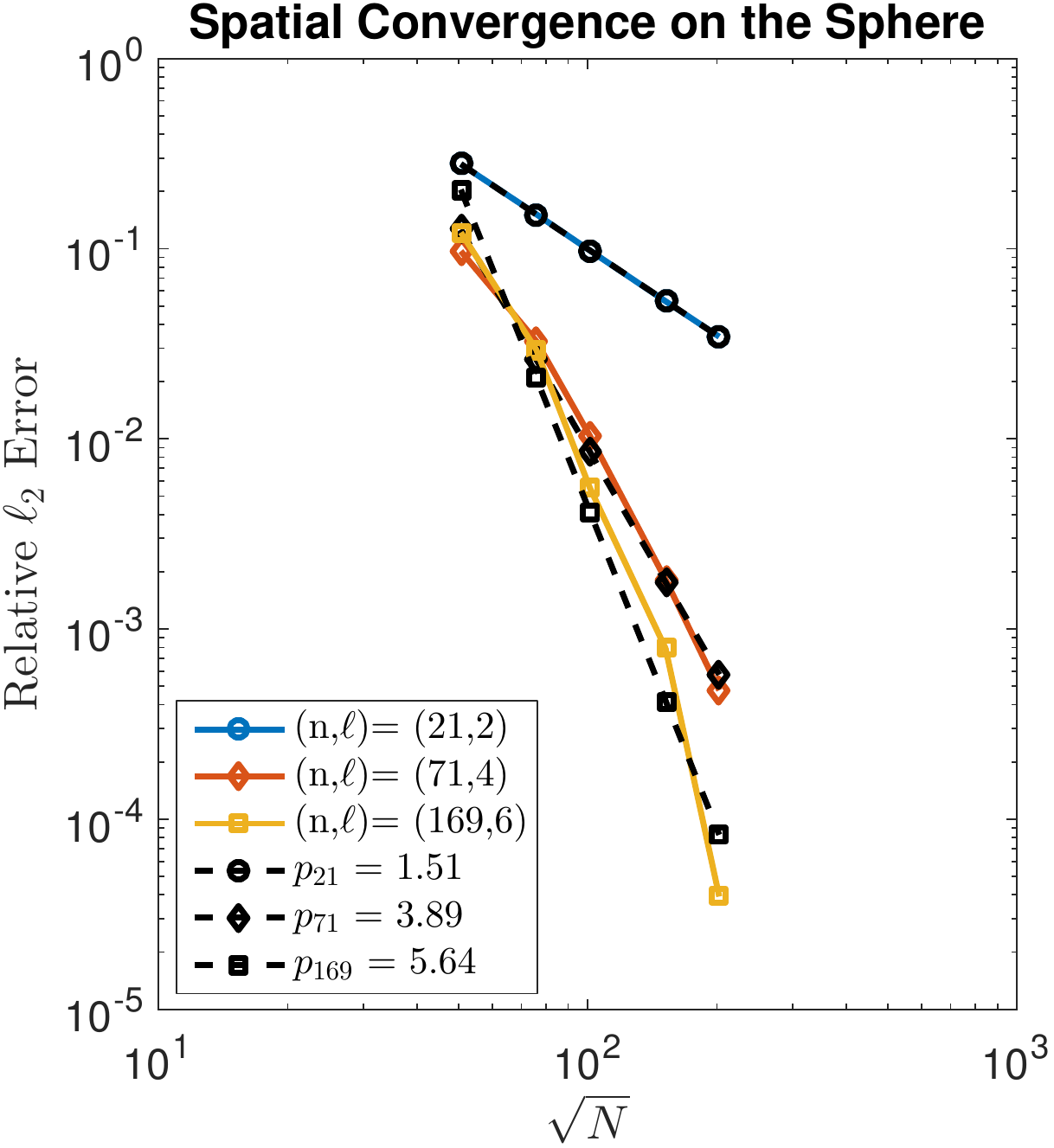}
\includegraphics[scale=0.6]{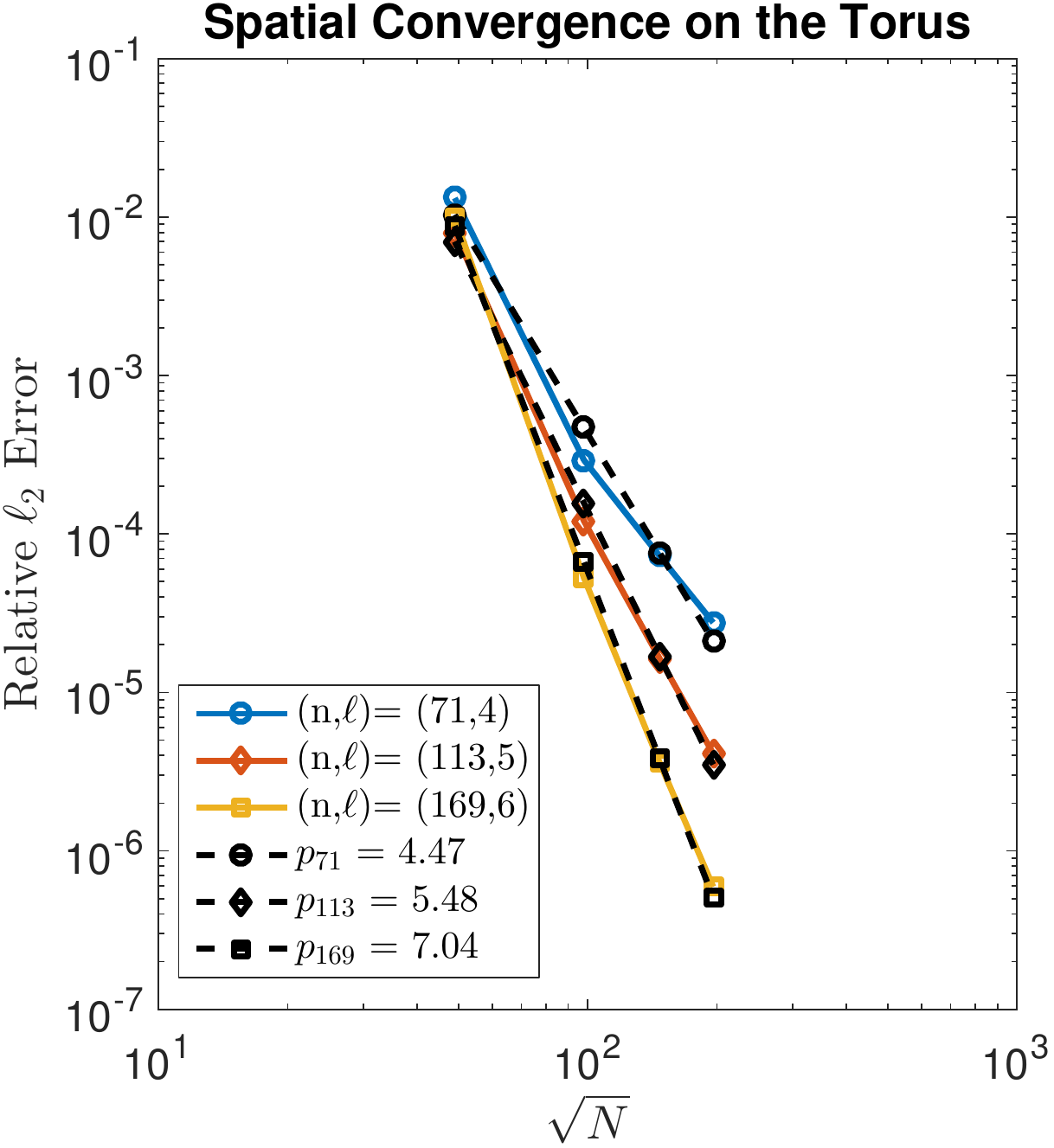}
\caption{Convergence on the sphere for the surface advection equation (left) and torus for the forced surface diffusion equation (right). The figure shows relative $\ell_2$ error as a function of $\sqrt{N}$ for different values of stencil $n$ and polynomial degree $\ell$.}
\label{fig:both}
\end{figure}

\subsection{Advection on the sphere $\mathbb{S}^2$}

In this test, we solve the surface advection equation on the sphere; in conservative form, this is given by:
\begin{align}
\frac{\partial c}{\partial t} + \nabla_{\mathbb{M}} \cdot \lf(\vu c\rt) = 0,
\label{eq:adv_sph}
\end{align}
where $c(\vx,t)$ is some scalar field being advected on the surface $\mathbb{M}$ by the velocity field $\vu(\vx,t)$. Of course, if $\nabla_{\mathbb{M}} \cdot \vu = 0$, this equation can be simplified to $\frac{\partial c}{\partial t} + \vu \cdot \nabla_{\mathbb{M}}c = 0$. However, we opt to discretize \eqref{eq:adv_sph} directly as this appears to produce lower errors than the alternative approach. Unfortunately, the RBF-FD differentiation matrices corresponding to the components of the surface divergence operator $\nabla_{\mathbb{M}}$ typically contain spurious eigenvalues in their spectra, often with (small) positive real parts. This can cause instabilities especially when an explicit time discretization is used. The remedy for this issue is to add a small amount of artificial hyperviscosity for stabilization~\cite{FoL11,FlyerLehto2012}. This transforms \eqref{eq:adv_sph} to
\begin{align}
\frac{\partial c}{\partial t} + \nabla_{\mathbb{M}} \cdot \lf(\vu c\rt) = \gamma \Delta_{\mathbb{M}}^k,
\label{eq:adv_sph_stab}
\end{align}
where $\gamma \in \mathbb{R}$ and $k \in \mathbb{N}$ must be tuned. Typically, $\gamma$ is a small number that goes to zero as $N \to \infty$, and $k$ is gently increased as the order of the method is increased~\cite{FlyerLehto2012}. For simplicity, we compute the discrete surface hyperviscosity operator by simply computing $L^k$, where $L$ is the discrete surface Laplacian. 

While formulas for $\gamma$ can be found in the literature for Euclidean domains in the context of PHS-based RBF-FD~\cite{FlyerNS,BarnettPHS}, there are no such formulas for the PHS-based RBF-FD on the sphere (to the best of our knowledge). In this work, we use the following formulas for $\gamma$ and $k$:
\begin{align}
\gamma &= (-1)^{1+k} 2^{2-2k} \lf(\sqrt{N}\rt)^{2-2k} \Lambda_{\rm max} \|\vu\|_{\rm max},\\
k &= \lf \lfloor \ln n \rt \rfloor,
\end{align}
where $\Lambda_{\rm max}$ is the \emph{real part} of the eigenvalue with largest real part of the sparse differentiation matrices $G^x,G^y$, and $G^z$. $\|\vu\|_{\rm max}$ is shorthand for the maximum of the pointwise $\ell_2$ norms of the velocities evaluated on the node set.  A derivation of this formula is beyond the scope of this paper, but will be detailed in future work. For this article, we always have $\|\vu\|_{\rm max} = 1$. $\Lambda_{\rm max}$ is estimated by Matlab calls to \rm{eigs(.,1,`LR')}, with a very loose tolerance of $8e-2$; this estimation is a preprocessing step for a given node set and stencil size. Our chosen test problem is the deformational flow test case from~\cite{NairLauritzen2010}. The components of the velocity field (in spherical coordinates) are
\begin{align}
u(\phi_1,\phi_2,t) &= \frac{10}{T} \cos\lf(\frac{\pi t}{T}\rt) \sin^2\lf(\phi_1 - \frac{2\pi t}{T} \rt)\sin\lf(2 \phi_2\rt) + \frac{2\pi}{T} \cos\lf(\phi_2\rt), \\
v(\phi_1,\phi_2,t) &= \frac{10}{T} \cos\lf(\frac{\pi t}{T}\rt) \sin\lf(2\phi_1 - \frac{2\pi t}{T} \rt)\cos\lf(\phi_2\rt),
\end{align}
where $-\pi \leq \phi_1 \leq \pi$, $-\pi/2 \leq \phi_2 \leq \pi/2$, and $T=5$. The flow field deforms the initial condition up to time $t=2.5$ and then reverses to return the solution to its initial position at $t=5$, which serves as the final time for the simulation. A simple change of basis is used to convert the velocity field into Cartesian coordinates. To test the convergence behavior of RBF-LOI under refinement, we use a smooth initial condition in the form of two Gaussian bells, given by
\begin{align}
c(\vx,0) = 0.95\lf(e^{-5\|\vx - {\bf p}_1\|_2^2} + e^{-5\|\vx - {\bf p}_2\|_2^2} \rt),
\end{align}
where ${\bf p}_1 = \lf(\sqrt{3}/2,1/2,0\rt)$ and ${\bf p}_2 = \lf(\sqrt{3}/2,-1/2,0\rt)$. Following~\cite{Aiton2011}, we use a time-step of $\Delta t = \frac{5}{2400}$ for this test. The time-stepping is done using the classical fourth-order explicit Runge Kutta method (RK4). The node sets were standard icosahedral nodes from the Spherepts package~\cite{spherepts}. The results are shown in Figure \ref{fig:both}(left). Figure \ref{fig:both} (left) shows that though one expects a convergence rate of $\xi = \ell$, we appear to obtain slightly lower rates. The LOI tolerance parameter $\tau$ was set to $\tau = 1e-2$ for $\ell = 2$, $\tau = 1e-3$ for $\ell = 3$, and $\tau = 1e-4$ for $\ell = 4$. In~\cite{FlyerLehto2012}, a shape parameter was carefully tuned to avoid stagnation errors, and the parameters $\gamma$ and $k$ were numerically computed so as to avoid instabilities. In contrast, our approach only involves setting $\tau$ and estimating $\Lambda_{\rm max}$ very crudely (which is done rapidly). No shape parameters or tuning were required to obtain stability, and no extended precision arithmetic was required. The goal of this article is simply demonstrate the feasibility of RBF-LOI for PDEs on surfaces, and we hence defer a deeper investigation of RBF-LOI for advection on the sphere to future work. We note that while the stencil sizes in our work are larger than those used in~\cite{FlyerLehto2012}, the increase in computational cost is largely ameliorated by the use of the overlapped RBF-FD method.

\subsection{Diffusion on a torus $\mathbb{T}$}

We consider the torus from~\cite{SWFKJSC2014} given by
\begin{align}
\mathbb{T} = \left\{\vX = (x,y,z) \in \mathbb{R}^3 \;\left| \left(1 - \sqrt{x^2 + y^2}\right)^2 + z^2 - \frac{1}{9} = 0 \right. \right\}.
\end{align}
Our goal is to solve the forced diffusion equation given by 
\begin{align}
\frac{\partial c}{\partial t} &= \Delta_{\mathbb{M}} c + f.
\end{align}
In all cases, we use the method of manufactured solutions, i.e., we prescribe a solution $c(\vx,t)$ and calculate the forcing term $f(\vx,t)$ that makes the solution hold. We use the BDF4 time-stepping scheme~\cite{Ascher97} for advancing the solution in time. This time-stepping scheme is fully implicit, and requires the solution of a sparse linear system every time-step. We set the time-step to $\Delta t = 10^{-3}$ for this test, and use Matlab's built-in sparse direct solver to solve the sparse linear systems obtained from overlapped RBF-FD. The manufactured solution in this case is given by
\begin{align}
c(t,\phi,\lambda) = e^{-5t} \sum\limits_{k=1}^{23}e^{-81(1-\cos(\lambda-\lambda_k))-9(1-\cos(\phi - \phi_k))},
\end{align}
where the intrinsic coordinates $-\pi\leq\phi,\lambda\leq\pi$ parameterize the torus $\mathbb{T}$ in the usual way~\cite{SWFKJSC2014}. The solution is $C^{\infty}(\mathbb{T})$. The LOI tolerance was fixed at $\tau = 10^{-3}$ for $\ell = 4,5$, and decreased to $\tau = 10^{-4}$ for $\ell=6$. The results are shown in Figure \ref{fig:both} (right). We obtain a convergence rate of approximately $\ell+1$ on the torus. Once again, no extended-precision arithmetic was used; all calculations are in double precision, with the cost of forming the larger differentiation matrices being almost completely ameliorated by the large speedup obtained from overlapped RBF-FD. Similar results were obtained for forced diffusion on the sphere as well (not shown).

\revone{
\subsection{Cost vs accuracy}
\begin{figure}[h!]
\includegraphics[scale=0.45]{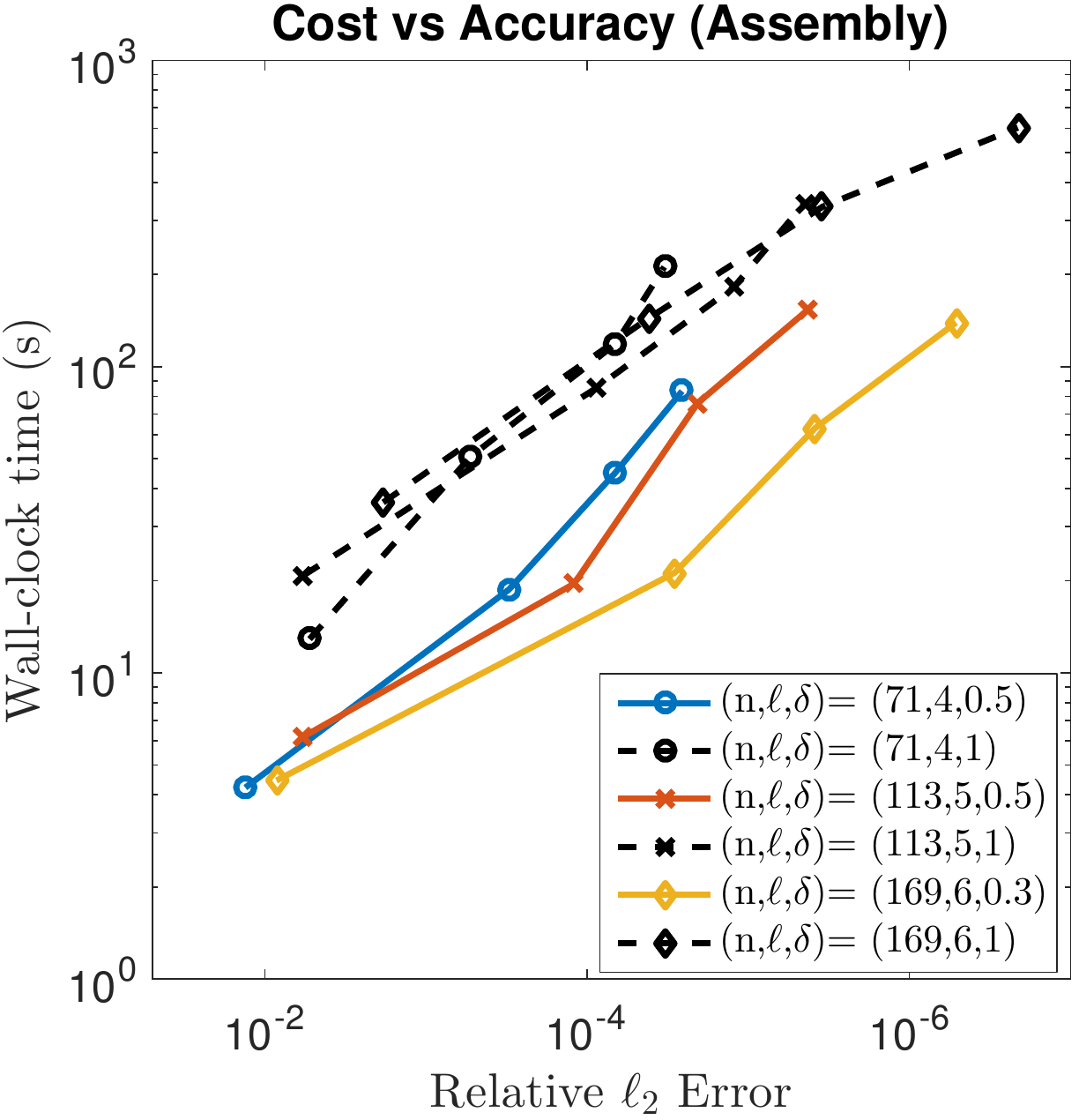}
\includegraphics[scale=0.45]{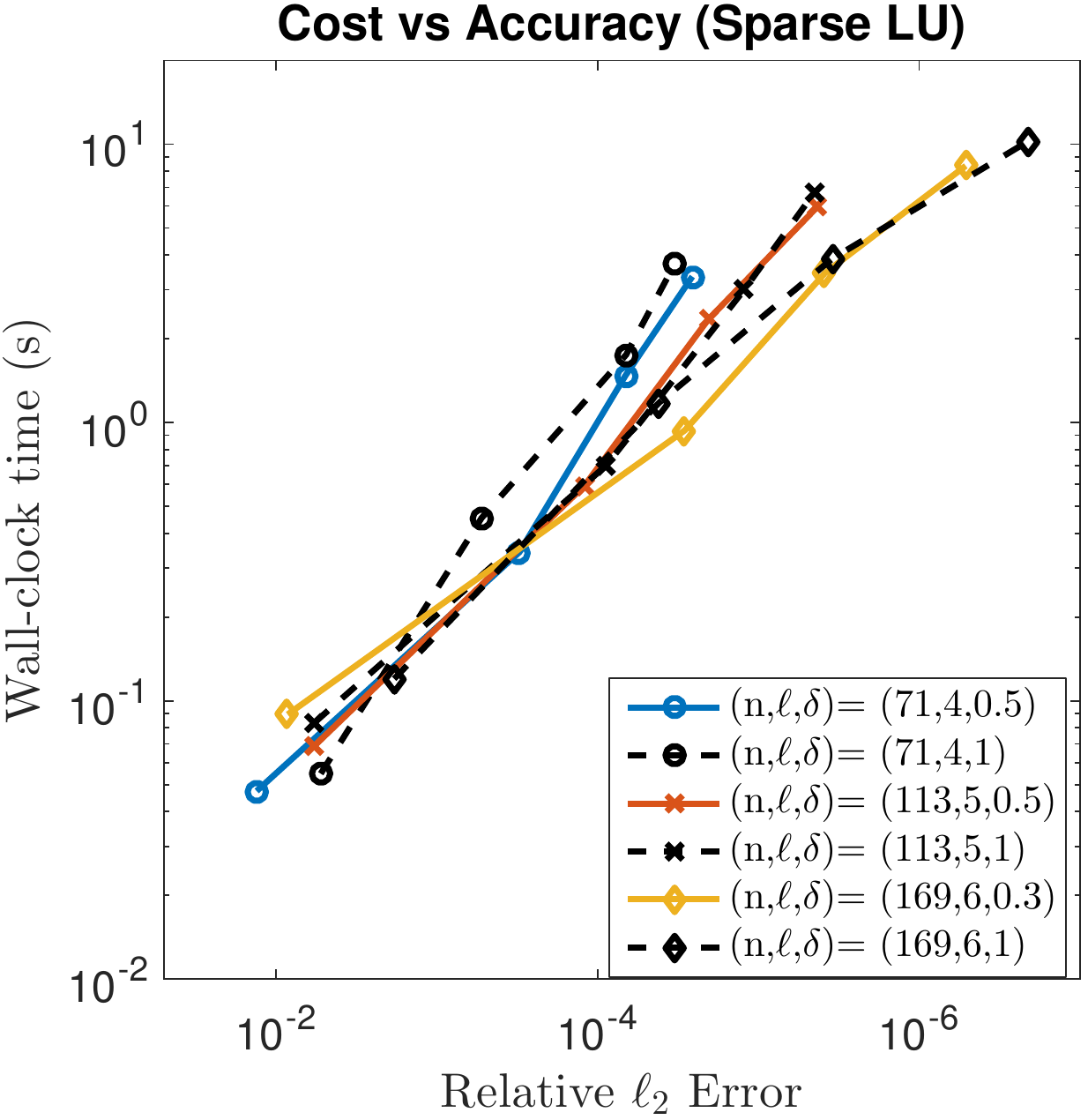}
\includegraphics[scale=0.45]{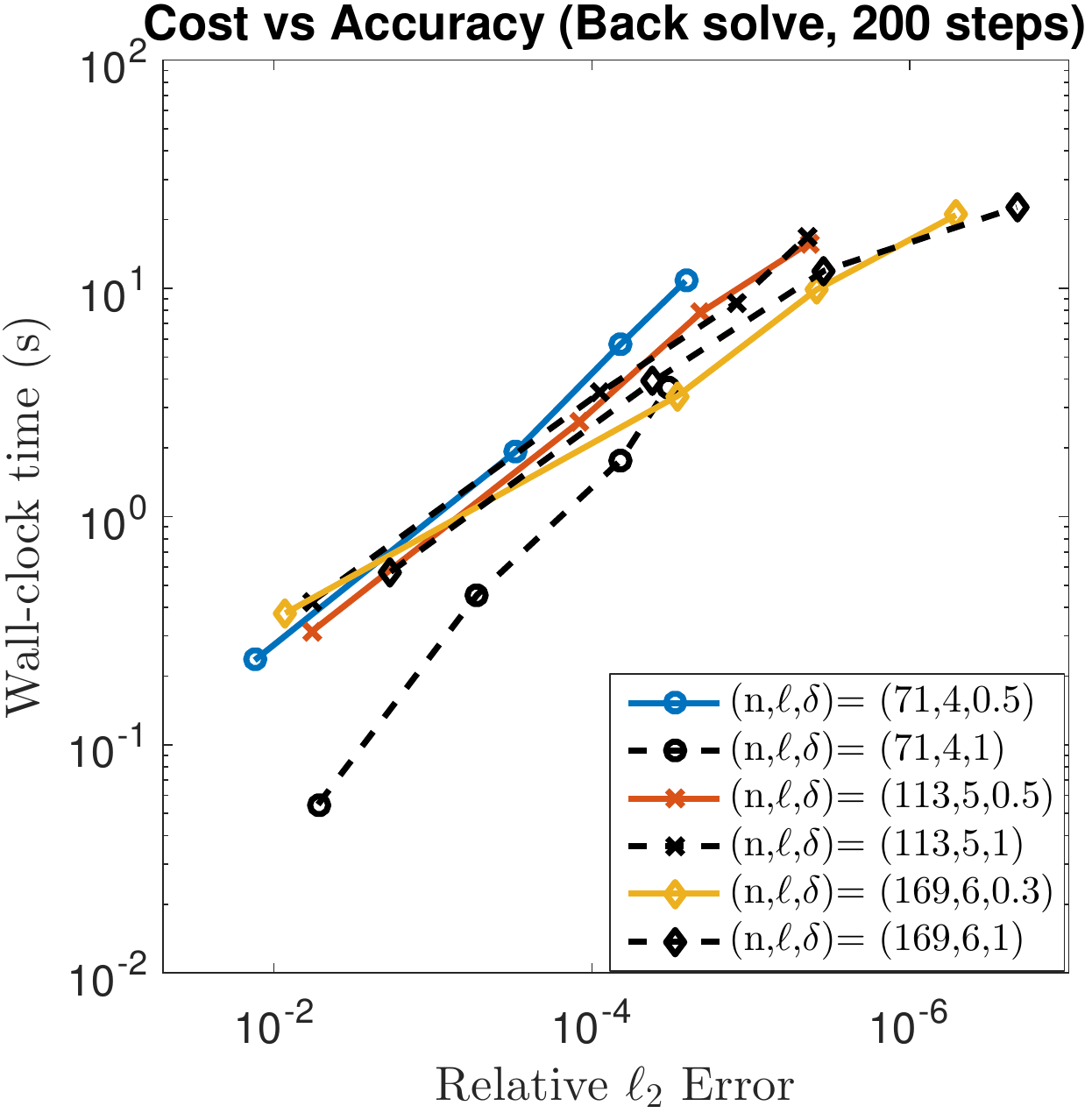}
\caption{Cost versus accuracy on the torus for the forced surface diffusion equation. The figures show wall-clock time as a function of relative $\ell_2$ error for different values of stencil $n$, polynomial degree $\ell$, and overlap parameter $\delta$. $\delta = 1$ corresponds to the standard RBF-FD method (dashed lines). The figures compare costs for the assembly stage (left), sparse LU factorization stage (middle), and the actual back-solve stage (right). In all cases, the $\ell_2$ error is a function of $\sqrt{N}$, which is increasing left to right.}
\label{fig:cva}
\end{figure}
To better understand the costs involved in both overlapped RBF-FD and its use within the RBF-LOI algorithm, we now study computational cost (measured in wall-clock time) as a function of accuracy for different values of $\xi$ (and therefore $n$ and $\ell$). Since the overlapped RBF-FD method is used only to speed up the assembly of the differentiation matrices, we present separate cost versus accuracy results for the assembly stage, the sparse LU factorization stage (preprocessing for solving the linear system), and the actual solution of the PDE using back solves, all with and without the use of overlapped RBF-FD. To ensure that the cost measured reflects matrix inversion rather than just matrix multiplies, we focus on the case of diffusion on the torus, which requires a sparse linear solve every step. The results are shown in Figure \ref{fig:cva}. 

Figure \ref{fig:cva} (left) shows that the overlapped RBF-FD method results in significant cost savings at the assembly stage for comparable accuracy; the method lags behind the standard RBF-FD method ($\delta = 1$) in accuracy only on the finest node set for $\ell=6$. Figure \ref{fig:cva} (middle) shows that the LU factorization costs are all comparable across different values of $\delta$, with the costs appearing to increase slowly with $n$ and $\ell$. Finally, Figure \ref{fig:cva} (right) shows that the back solve costs are also mostly comparable across all methods, with $\ell=4$ being the exception; in this case, $\delta=1$ appears to be cheaper, possibly due to slightly different matrix structure. It is important to note that if one sums up costs across the subfigures of Figure \ref{fig:cva} for $\delta = 1$, the assembly cost dominates all other costs. In contrast, for the overlapped RBF-FD method, the assembly cost is comparable to the back-solve cost. This feature would be beneficial when solving a problem on a moving domain.
}

%% file: Applications.tex
\section{Applications}
\label{sec:appl}

Having validated our the RBF-LOI method on standard test cases, we now turn our attention to some applications. Our goal here is demonstrate that the RBF-LOI method is stable on different point cloud surfaces \emph{and} PDEs that are more complicated than the forced diffusion equation. To that end, we test on three manifolds of increasing genus:
\begin{enumerate}
\item The red blood cell (genus 0), a parametric surface with node sets and normals obtained using the techniques outlined in~\cite{SFKSISC2017};
\item Dupin's cyclide (genus 1), an implicit surface with node sets and normals obtained using Meshlab~\cite{Meshlab}; and
\item The double torus (genus 2), another implicit surface with node sets and normals again obtained using Meshlab.
\end{enumerate}
Unfortunately, unlike in~\cite{SWFKJSC2014,LSWSISC2017}, we were unable to find stable parameters for point cloud models of more complicated manifolds such as frogs and bunnies. It is likely that such surfaces would require an adaptive tolerance selection for the LOI procedure, which we leave for future work. The problem of advection on arbitrary surfaces also requires a very careful derivation of hyperviscosity parameters $\gamma$ and $k$, which is likewise beyond the scope of this article. We instead focus on biologically-motivated reaction-diffusion models involving nonlinear terms. We believe that these applications serve as a convincing demonstration of the simplicity and effectiveness of the RBF-LOI method. In all cases, we use the following parameters for RBF-LOI: $\ell=4$, $m=2\ell+1$, tolerance of $\tau = 10^{-3}$. All simulations used the SBDF2 method for time-stepping, with corresponding linear systems being solved by the built-in Matlab sparse direct solver.

\subsection{Cahn-Hilliard on a Red Blood Cell}
\begin{figure}[h!]
\subfloat[]
{
\includegraphics[scale=0.6]{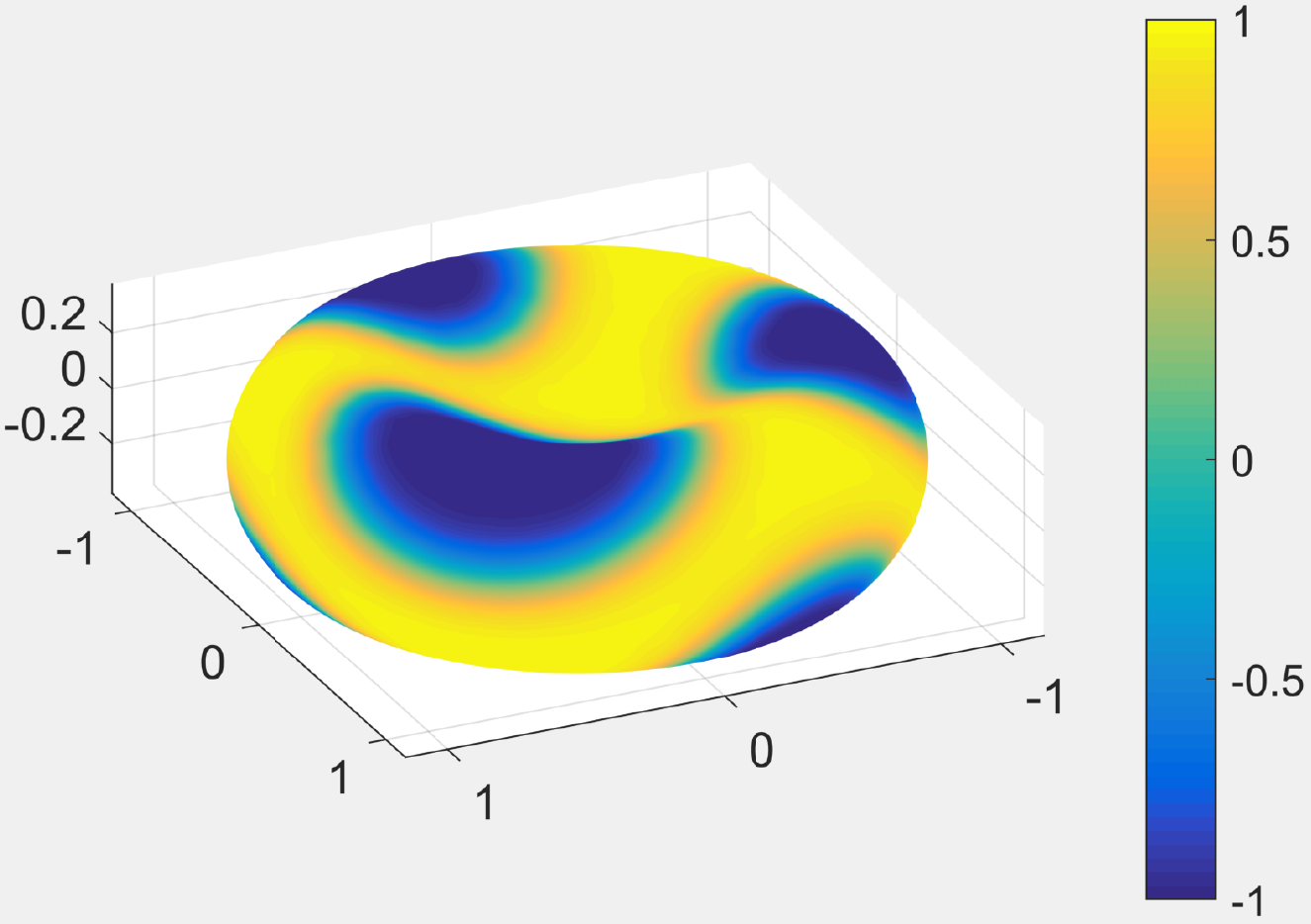}
}
\subfloat[]
{
\includegraphics[scale=0.7]{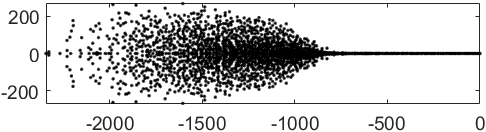}
}
\caption{Solution of the Cahn-Hilliard equation on the red blood cell at time $T=2$ (left) and eigenvalues of the Laplacian (right).}
\label{fig:ch}
\end{figure}

Our first application is the simulation of the Cahn-Hilliard equation on an idealized red blood cell~\cite{FuselierWright:2013}. The Cahn-Hilliard equation is a nonlinear PDE governing phase separation, and is given by
\begin{align}
\frac{\partial c}{\partial t} &= \nu \Delta_{\mathbb{M}}c^3 - \nu\Delta_{\mathbb{M}}c - \nu\gamma\Delta^2_{\mathbb{M}} c,
\end{align}
where $\Delta^2_{\mathbb{M}}$ is the surface bilaplacian. The solutions $c=1$ and $c=-1$ both constitute critical points of this reaction-diffusion system, and any initial condition will be separated over time into these two phases. We simulate the above PDE on the red blood cell to time $t=2$ using a random initial condition. To approximate the surface bilaplacian, we first form the discrete surface Laplacian $L$, then simply compute the discrete surface bilaplacian $B$ as $B=L.L$. This has the effect of increasing the fill-in of $B$ when compared to $L$, but our goal here is to simply demonstrate effectiveness. We use $\gamma = 0.006$ and $\nu = 0.5$, and a time-step of $\Delta t= 10^{-4}$. This small step is primarily due to the stiff nonlinear term $\nu \Delta_{\mathbb{M}}c^3$ being stepped explicitly in time. The results for $N=2553$ nodes are shown in Figure \ref{fig:ch}(a), and the spectrum of the discrete surface Laplacian is shown in Figure \ref{fig:ch}(b). Clearly, our solutions are exhibiting the correct qualitative behavior, and the spectrum of $L$ contains no spurious eigenvalues. To obtain the same behavior without LOI, careful tuning of the shape parameter on a stencil-by-stencil basis was needed in~\cite{SWFKJSC2014}, and a stencil selection algorithm was needed in~\cite{LSWSISC2017}.

\subsection{Fitzhugh-Nagumo waves on Dupin's cyclide}
\begin{figure}[h!]
\subfloat[]
{
\includegraphics[scale=0.4]{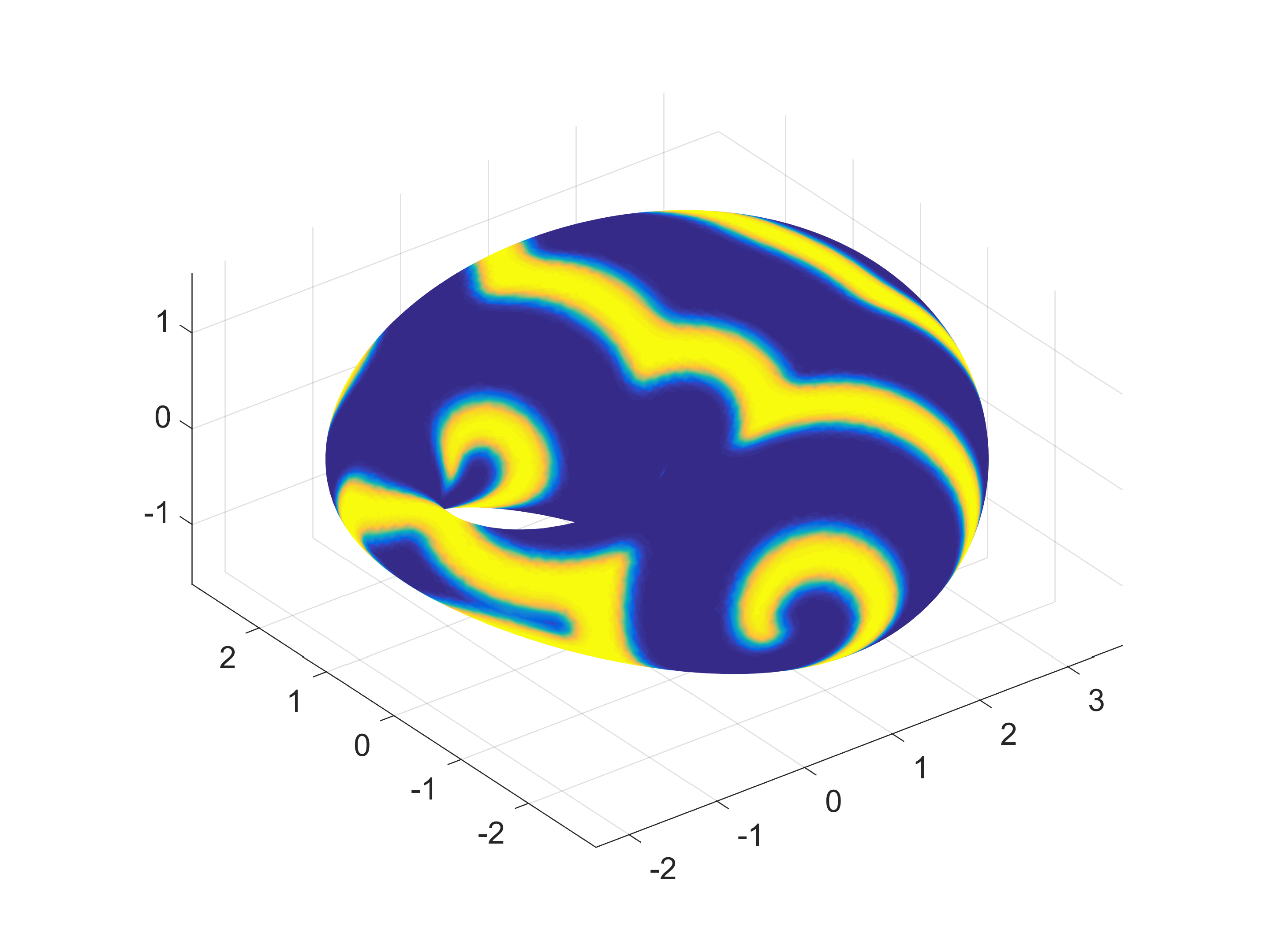}
}
\subfloat[]
{
\includegraphics[scale=0.7]{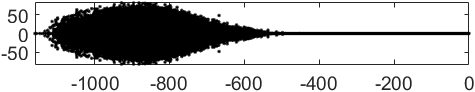}
}
\caption{Solution of the Fitzhugh-Nagumo equations on Dupin's cyclide at time $T=800$ (left) and eigenvalues of the Laplacian (right).}
\label{fig:dupin}
\end{figure}
Following~\cite{FuselierWright:2013}, we simulate the Fitzhugh-Nagumo reaction-diffusion system on Dupin's cyclide. The reaction-diffusion model is given by:
\begin{align}
\frac{\partial c_1}{\partial t} &= \delta_1 \Delta_{\mathbb{M}} c_1 + \frac{1}{0.02} c_1 \lf(1-c_1\rt)\lf(c_1 - \frac{c_2 + 0.02}{0.75}\rt), \\
\frac{\partial c_2}{\partial t} &= \delta_1 \Delta_{\mathbb{M}} c_2 + c_1 - c_2,
\end{align}
where $c_1$ and $c_2$ are typically viewed as chemical concentrations or densities corresponding to a membrane potential and a current, respectively. The above system is a simple model for the dynamics of excitable media, and is often viewed as a simplification of the Hodgkin-Huxley model for the dynamics of neurons~\cite{fitzhugh1961,nagumo1962}. Our initial condition on Dupin's cyclide is given by $c_1 = \frac{1}{2}\lf(1 + \tanh(5x + y)\rt)$ and $c_2 = \frac{1}{2}\lf(1 - \tanh(10z) \rt)$, where $\bs{x} = (x,y,z)$. The node sets on the cyclide are the same as those used in~\cite{FuselierWright:2013}; we use $N= 11884$ of these nodes. The results of the simulation with the SBDF2 method at time $t=100$ are shown in Figure \ref{fig:dupin}(a), and the spectrum of the discrete Laplacian is shown in Figure \ref{fig:dupin}(b). As expected, the simulation results in spiral waves that scroll over the manifold. 

\subsection{Turing spots on the double torus}
\begin{figure}[h!]
\subfloat[]
{
\includegraphics[scale=0.6]{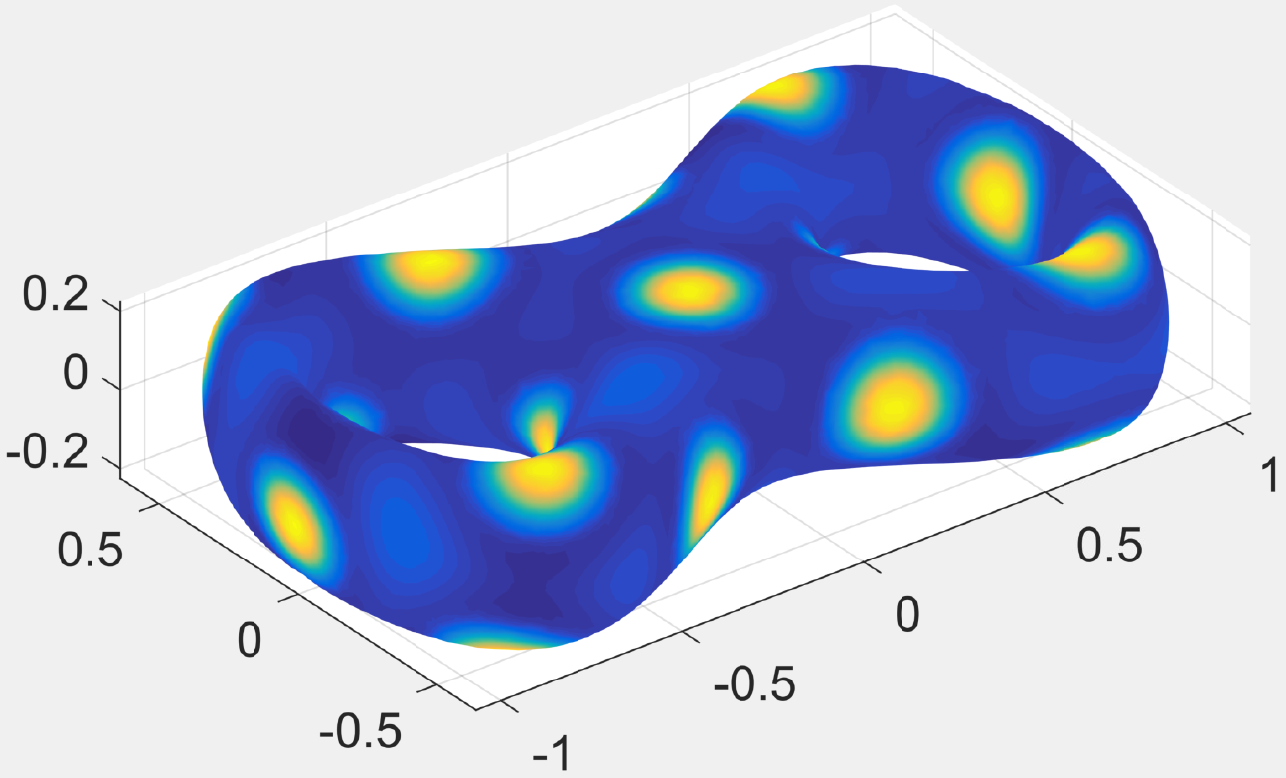}
}
\subfloat[]
{
\includegraphics[scale=0.7]{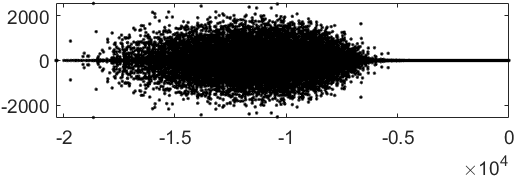}
}
\caption{Solution of the Turing equations on the double torus at time $T=800$ (left) and eigenvalues of the Laplacian (right).}
\label{fig:turing}
\end{figure}
Our final application involves solving another coupled reaction-diffusion system on the double torus
\begin{align}
\mathbb{T}^2 = \{ \vx = (x,y,z) \in \mathbb{R}^3 \mid (x^2(1-x^2) - y^2)^2 + 0.5z^2 = \frac{1}{40}\},
\end{align}
which a genus-2 surface obtained as the \emph{join} of two genus-1 tori. On this surface, we simulate the Turing system given by
\begin{align}
\frac{\partial c_1}{\partial t} &= \delta_1 \Delta_{\mathbb{M}}c_1 + \alpha c_1 \lf(1 - \tau_1 c_2^2\rt) + c_2 \lf(1-\tau_2 c_1\rt), \\
\frac{\partial c_2}{\partial t} &= \delta_2 \Delta_{\mathbb{M}}c_2 + \beta c_2 \lf(1 + \frac{\alpha \tau_1}{\beta} c_1 c_2\rt) + c_1 \lf(\gamma_1 + \tau_2c_2 \rt),
\end{align}
where we use the parameters $\delta_1 = 0.0011$, $\delta_2 = 0.0021$, $\tau_1 = 0.02$, $\tau_2 = 0.2$, $\alpha = 0.899$, $\beta = -0.91$, and $\gamma_1 = -\alpha$. We use a time-step of $\Delta t = 0.01$ and simulate to a final time of $t= 800$ on $N=12100$ nodes. The results are shown in Figure \ref{fig:turing}a, and the spectrum of the discrete Laplacian $L$ is shown in Figure \ref{fig:turing}b. Figure \ref{fig:turing}a shows that spot patterns have formed on the double torus despite the relatively coarse spatial discretization.

%% file: Discussion.tex
\section{Summary and Future Work}
\label{sec:summary}

We have proposed and demonstrated numerical solutions to PDEs on manifolds via RBF-LOI: a polynomially-augmented RBF-FD procedure. The major novel contribution of our work has been demonstration that a well-chosen polynomial basis (here, the Least Orthogonal Interpolant) along with efficient stencil overlap techniques for RBF-FD approximation can yield a stable, robust, efficient, and accurate PDE solver on manifolds. Our algorithm relies on specification of only two global parameters: an overlap parameter (which affects cost of local discrete operator construction), and a LOI tolerance parameter (which affects stability of polynomial approximations). 

To tackle more general point clouds, we would like the LOI tolerance parameter to be automatically tunable on a per-stencil basis. Ongoing work revolves around devising an automated approach for this parameter, and application of the RBF-LOI to PDE solutions on more intricate manifolds. Future work would involve rigorously deriving a hyperviscosity formulation for PDEs on arbitrary orientable manifolds to facilitate the solution of hyperbolic problems with RBF-LOI.

%% file: article.bbl
\begin{thebibliography}{}

\bibitem[\protect\astroncite{Aiton}{2014}]{Aiton2011}
Aiton, K.~A. (2014).
\newblock {A Radial Basis Function Partition of Unity Method for Transport on
  the Sphere}.
\newblock Master's thesis, Boise State University, USA.

\bibitem[\protect\astroncite{Ascher et~al.}{1997}]{Ascher97}
Ascher, U.~M., Ruuth, S.~J., and Wetton, B. T.~R. (1997).
\newblock {Implicit-Explicit Methods For Time-Dependent PDEs}.
\newblock {\em SIAM J. Numer. Anal}, 32:797--823.

\bibitem[\protect\astroncite{Barnett}{2015}]{BarnettPHS}
Barnett, G.~A. (2015).
\newblock {\em A Robust RBF-FD Formulation based on Polyharmonic Splines and
  Polynomials}.
\newblock PhD thesis, University of Colorado Boulder.

\bibitem[\protect\astroncite{Bayona et~al.}{2017}]{FlyerElliptic}
Bayona, V., Flyer, N., Fornberg, B., and Barnett, G.~A. (2017).
\newblock On the role of polynomials in rbf-fd approximations: Ii. numerical
  solution of elliptic pdes.
\newblock {\em Journal of Computational Physics}, 332(Supplement C):257 -- 273.

\bibitem[\protect\astroncite{Bayona et~al.}{2010}]{Bayona2010}
Bayona, V., Moscoso, M., Carretero, M., and Kindelan, M. (2010).
\newblock {RBF-FD} formulas and convergence properties.
\newblock {\em J. Comput. Phys.}, 229(22):8281--8295.

\bibitem[\protect\astroncite{Behrens and Iske}{2002}]{iske2002}
Behrens, J. and Iske, A. (2002).
\newblock {Grid-free adaptive semi-{L}agrangian advection using radial basis
  functions}.
\newblock {\em Comput. Math. Appl.}, 43(3):319--327.

\bibitem[\protect\astroncite{Boor and Ron}{1992}]{boor_computational_1992}
Boor, C.~D. and Ron, A. (1992).
\newblock Computational {Aspects} of {Polynomial} {Interpolation} in {Several}
  {Variables}.
\newblock {\em Mathematics of Computation}, 58(198):705--727.

\bibitem[\protect\astroncite{Cignoni et~al.}{2008}]{Meshlab}
Cignoni, P., Callieri, M., Corsini, M., Dellepiane, M., Ganovelli, F., and
  Ranzuglia, G. (2008).
\newblock {MeshLab: an Open-Source Mesh Processing Tool}.
\newblock In Scarano, V., Chiara, R.~D., and Erra, U., editors, {\em
  Eurographics Italian Chapter Conference}. The Eurographics Association.

\bibitem[\protect\astroncite{Davydov and Oanh}{2011}]{Davydov2011}
Davydov, O. and Oanh, D.~T. (2011).
\newblock Adaptive meshless centres and {RBF} stencils for poisson equation.
\newblock {\em J. Comput. Phys.}, 230(2):287--304.

\bibitem[\protect\astroncite{Davydov and Schaback}{2017}]{DavydovSchaback2017}
Davydov, O. and Schaback, R. (2017).
\newblock Optimal stencils in sobolev spaces.
\newblock Submitted.

\bibitem[\protect\astroncite{Fasshauer}{2007}]{Fasshauer:2007}
Fasshauer, G.~E. (2007).
\newblock {\em {Meshfree Approximation Methods with {MATLAB}}}.
\newblock Interdisciplinary Mathematical Sciences - Vol. 6. World Scientific
  Publishers, Singapore.

\bibitem[\protect\astroncite{Fasshauer and McCourt}{2012}]{FaMC12}
Fasshauer, G.~E. and McCourt, M.~J. (2012).
\newblock Stable evaluation of {G}aussian radial basis function interpolants.
\newblock {\em SIAM J. Sci. Comput.}, 34:A737----A762.

\bibitem[\protect\astroncite{FitzHugh}{1961}]{fitzhugh1961}
FitzHugh, R. (1961).
\newblock Impulses and physiological states in theoretical models of nerve
  membrane.
\newblock {\em Biophysical journal}, 1(6):445--466.

\bibitem[\protect\astroncite{Flyer et~al.}{2016a}]{FlyerNS}
Flyer, N., Barnett, G.~A., and Wicker, L.~J. (2016a).
\newblock Enhancing finite differences with radial basis functions: Experiments
  on the {N}avier-{S}tokes equations.
\newblock {\em J. Comput. Phys.}, 316:39--62.

\bibitem[\protect\astroncite{Flyer et~al.}{2016b}]{FlyerPHS}
Flyer, N., Fornberg, B., Bayona, V., and Barnett, G.~A. (2016b).
\newblock {On the role of polynomials in RBF-FD approximations: I.
  Interpolation and accuracy}.
\newblock {\em J. Comput. Phys.}, 321:21--38.

\bibitem[\protect\astroncite{Flyer et~al.}{2012}]{FlyerLehto2012}
Flyer, N., Lehto, E., Blaise, S., Wright, G.~B., and St-Cyr, A. (2012).
\newblock {A guide to {RBF}-generated finite differences for nonlinear
  transport: shallow water simulations on a sphere}.
\newblock {\em J. Comput. Phys.}, 231:4078--4095.

\bibitem[\protect\astroncite{Flyer and Wright}{2007}]{FlyerWright:2007}
Flyer, N. and Wright, G.~B. (2007).
\newblock Transport schemes on a sphere using radial basis functions.
\newblock {\em J. Comput. Phys.}, 226:1059--1084.

\bibitem[\protect\astroncite{Flyer and Wright}{2009}]{FlyerWright:2009}
Flyer, N. and Wright, G.~B. (2009).
\newblock {A radial basis function method for the shallow water equations on a
  sphere}.
\newblock {\em Proc. Roy. Soc. A}, 465:1949--1976.

\bibitem[\protect\astroncite{Fornberg et~al.}{2011}]{FLF}
Fornberg, B., Larsson, E., and Flyer, N. (2011).
\newblock Stable computations with {G}aussian radial basis functions.
\newblock {\em SIAM J. Sci. Comput.}, 33(2):869--892.

\bibitem[\protect\astroncite{Fornberg and Lehto}{2011}]{FoL11}
Fornberg, B. and Lehto, E. (2011).
\newblock {Stabilization of {RBF}-generated finite difference methods for
  convective {PDE}s}.
\newblock {\em J. Comput. Phys.}, 230:2270--2285.

\bibitem[\protect\astroncite{Fornberg et~al.}{2013}]{FoLePo13}
Fornberg, B., Lehto, E., and Powell, C. (2013).
\newblock Stable calculation of {G}aussian-based {RBF-FD} stencils.
\newblock {\em Comput. Math. Appl.}, 65:627--637.

\bibitem[\protect\astroncite{Fornberg and Piret}{2007}]{FornbergPiret:2007}
Fornberg, B. and Piret, C. (2007).
\newblock {A stable algorithm for flat radial basis functions on a sphere}.
\newblock {\em SIAM J. Sci. Comput.}, 30:60--80.

\bibitem[\protect\astroncite{Fornberg and Wright}{2004}]{FoWr}
Fornberg, B. and Wright, G. (2004).
\newblock {Stable computation of multiquadric interpolants for all values of
  the shape parameter}.
\newblock {\em Comput. Math. Appl.}, 48:853--867.

\bibitem[\protect\astroncite{Fuselier and Wright}{2013}]{FuselierWright:2013}
Fuselier, E.~J. and Wright, G.~B. (2013).
\newblock {A high-order kernel method for diffusion and reaction-diffusion
  equations on surfaces}.
\newblock {\em J. Sci. Comput.}, 56(3):535--565.

\bibitem[\protect\astroncite{Lehto et~al.}{2017}]{LSWSISC2017}
Lehto, E., Shankar, V., and Wright, G.~B. (2017).
\newblock A radial basis function (rbf) compact finite difference (fd) scheme
  for reaction-diffusion equations on surfaces.
\newblock {\em SIAM Journal on Scientific Computing}, 39(5):A2129--A2151.

\bibitem[\protect\astroncite{Nagumo et~al.}{1962}]{nagumo1962}
Nagumo, J., Arimoto, S., and Yoshizawa, S. (1962).
\newblock An active pulse transmission line simulating nerve axon.
\newblock {\em Proceedings of the IRE}, 50(10):2061--2070.

\bibitem[\protect\astroncite{Nair and Lauritzen}{2010}]{NairLauritzen2010}
Nair, R.~D. and Lauritzen, P.~H. (2010).
\newblock A class of deformational flow test cases for linear transport
  problems on the sphere.
\newblock {\em J. Comput. Phys.}, 229(23):8868--8887.

\bibitem[\protect\astroncite{Narayan and Xiu}{2012}]{narayan_stochastic_2012}
Narayan, A. and Xiu, D. (2012).
\newblock Stochastic {Collocation} {Methods} on {Unstructured} {Grids} in
  {High} {Dimensions} via {Interpolation}.
\newblock {\em SIAM Journal on Scientific Computing}, 34(3):A1729--A1752.

\bibitem[\protect\astroncite{Piret}{2012}]{Piret2012}
Piret, C. (2012).
\newblock {The orthogonal gradients method: A radial basis functions method for
  solving partial differential equations on arbitrary surfaces}.
\newblock {\em J. Comput. Phys.}, 231(20):4662--4675.

\bibitem[\protect\astroncite{Piret and Dunn}{2016}]{Piret2016}
Piret, C. and Dunn, J. (2016).
\newblock Fast rbf ogr for solving pdes on arbitrary surfaces.
\newblock {\em AIP Conference Proceedings}, 1776(1).

\bibitem[\protect\astroncite{Reeger and Fornberg}{2016}]{ReegerFornbergQuad1}
Reeger, J.~A. and Fornberg, B. (2016).
\newblock Numerical quadrature over the surface of a sphere.
\newblock {\em Studies in Applied Mathematics}, 137(2):174--188.

\bibitem[\protect\astroncite{Reeger and Fornberg}{2018}]{ReegerFornbergQuad3}
Reeger, J.~A. and Fornberg, B. (2018).
\newblock Numerical quadrature over smooth surfaces with boundaries.
\newblock {\em Journal of Computational Physics}, 355(Supplement C):176 -- 190.

\bibitem[\protect\astroncite{Reeger et~al.}{2016}]{ReegerFornbergQuad2}
Reeger, J.~A., Fornberg, B., and Watts, M.~L. (2016).
\newblock Numerical quadrature over smooth, closed surfaces.
\newblock {\em Proceedings of the Royal Society of London A: Mathematical,
  Physical and Engineering Sciences}, 472(2194).

\bibitem[\protect\astroncite{Schaback}{2005}]{Schaback2005}
Schaback, R. (2005).
\newblock Multivariate interpolation by polynomials and radial basis functions.
\newblock {\em Constructive Approximation}, 21(3):293--317.

\bibitem[\protect\astroncite{Shankar}{2017}]{ShankarJCP2017}
Shankar, V. (2017).
\newblock The overlapped radial basis function-finite difference ({RBF-FD})
  method: A generalization of {RBF-FD}.
\newblock {\em J. Comput. Phys.}, 342:211--228.

\bibitem[\protect\astroncite{Shankar et~al.}{2017}]{SFKSISC2017}
Shankar, V., Kirby, R.~M., and Fogelson, A.~L. (2017).
\newblock Robust node generation for meshfree discretizations on irregular
  domains and surfaces.
\newblock Submitted.

\bibitem[\protect\astroncite{Shankar and Wright}{2018}]{SWJCP2018}
Shankar, V. and Wright, G.~B. (2018).
\newblock Mesh-free semi-lagrangian methods for transport on a sphere using
  radial basis functions.
\newblock {\em J. Comput. Phys.}, 366(C):170--190.

\bibitem[\protect\astroncite{Shankar et~al.}{2014}]{SWFKJSC2014}
Shankar, V., Wright, G.~B., Kirby, R.~M., and Fogelson, A.~L. (2014).
\newblock A radial basis function ({RBF})-finite difference ({FD}) method for
  diffusion and reaction--diffusion equations on surfaces.
\newblock {\em J. Sci. Comput.}, 63(3):745--768.

\bibitem[\protect\astroncite{Wendland}{2005}]{Wendland:2004}
Wendland, H. (2005).
\newblock {\em {Scattered data approximation}}, volume~17 of {\em Cambridge
  Monogr. Appl. Comput. Math.}
\newblock Cambridge University Press, Cambridge.

\bibitem[\protect\astroncite{Wright}{2018}]{spherepts}
Wright, G.~B. (2018).
\newblock {SpherePts}.
\newblock {https://github.com/gradywright/spherepts/}.

\bibitem[\protect\astroncite{Wright and Fornberg}{2006}]{Wright200699}
Wright, G.~B. and Fornberg, B. (2006).
\newblock {Scattered node compact finite difference-type formulas generated
  from radial basis functions}.
\newblock {\em J. Comput. Phys.}, 212(1):99--123.

\bibitem[\protect\astroncite{Wright and Fornberg}{2017}]{FoWr2016}
Wright, G.~B. and Fornberg, B. (2017).
\newblock Stable computations with flat radial basis functions using
  vector-valued rational approximations.
\newblock {\em J. Comput. Phys.}, 331:137 -- 156.

\end{thebibliography}
